\documentclass[a4paper,12pt]{article}

\usepackage{latexsym}
\usepackage{amssymb}
\usepackage{theorem}
\usepackage{amsmath}
\usepackage{amscd}
\usepackage{graphicx}
\usepackage{color}

\newtheorem{theorem}{Theorem}[section]

\newtheorem{example}[theorem]{Example}
\newtheorem{proposition}[theorem]{Proposition}
\newtheorem{remark}[theorem]{Remark}
\newtheorem{definition}[theorem]{Definition}

\newcommand{\demo}{\par\noindent{\it Proof. \/}\ }
\newcommand{\enD}{\hfill $\Box$\vspace{3truemm} \par}
\newcommand{\R}{\mathbb{R}}
\newcommand{\Z}{\mathbb{Z}}

\newcommand{\bn}{\mbox{\boldmath $n$}}
\newcommand{\bt}{\mbox{\boldmath $t$}}
\newcommand{\ba}{\mbox{\boldmath $a$}}
\newcommand{\bb}{\mbox{\boldmath $b$}}
\newcommand{\bc}{\mbox{\boldmath $c$}}
\newcommand{\be}{\mbox{\boldmath $e$}}
\newcommand{\bmu}{\mbox{\boldmath $\mu$}}
\newcommand{\bv}{\mbox{\boldmath $v$}}
\newcommand{\bw}{\mbox{\boldmath $w$}}

\begin{document}

\title{Bertrand types of regular curves and Bertrand framed curves in the Euclidean 3-space}

\author{Nozomi Nakatsuyama and Masatomo Takahashi}

\date{\today}

\maketitle

\begin{abstract}
A Bertrand (respectively, Mannheim) curve is a space curve whose principal normal line is the same as the principal normal (respectively, bi-normal) line of another curve. 
By definition, another curve is a parallel curve with respect to the direction of the principal normal vector. 
In this paper, we consider the other cases, that is, a space curve whose tangent (or, principal normal, bi-normal) line is the same as the tangent (or, principal normal, bi-normal) line of another curve, respectively. 
We say that a Bertrand type curve if there exists such another curve. 
We clarify that the existence conditions of Bertrand type curves in all cases. 
There are times when the Bertrand type curve does not exist. 
On the other hand, since the another curve may have singular points, we also consider curves with singular points. 
As smooth curves with singular points, it is useful to use the framed curves in the Euclidean space.
Then we define and investigate Bertrand framed curves. 
We also clarify that the existence conditions of the Bertrand framed curves in all cases.
\end{abstract}

\renewcommand{\thefootnote}{\fnsymbol{footnote}}

%%%%%%%%%%%%%%%% Section 1 %%%%%%%%%%%%%%%%%%
\section{Introduction}

Bertrand and Mannheim curves are classical objects in differential geometry (\cite{Aminov, Banchoff-Lovett, Berger-Gostiaux, Bertrand, doCarmo, HCIP, Izumiya-Takeuchi1, Kuhnel, Liu-Wang, Struik}). 
A Bertrand (respectively, Mannheim) curve is a space curve whose principal normal line is the same as the principal normal (respectively, bi-normal) line of another curve. 
By definition, another curve is a parallel curve with respect to the direction of the principal normal vector. 
In order to define the principal normal vector, the non-degenerate condition is needed.
In general, the parallel curve does not satisfy the non-degenerate condition. 
Even if regular cases, the existence conditions of the Bertrand and Mannheim curves seem to be missed the non-degenerate condition in some books and papers.  
Bertrand curves have been applied in computer-aided geometric design (cf. \cite{Papaioannou-Kiritsis}). 
In \cite{Honda-Takahashi-2020}, we investigated Bertrand and Mannheim curves of non-degenerate curves under the assumption that the torsion does not vanish. 
In this paper, we give existence conditions of Bertrand and Mannheim curves without this assumption  in \S 2. 
Since we have three lines, that is, the tangent, principal normal and bi-normal lines, it is natural to ask  how about the other cases. 
We consider the other cases, that is, a space curve whose tangent (or, principal normal, bi-normal) line is the same as the tangent (or, principal normal, bi-normal) line of another curve, respectively. 
We say that a Bertrand type curve if there exists such another curve. 
In \S 3, we clarify that the existence conditions of Bertrand type curves in all cases. 
As a consequence, there are times when the Bertrand type curve does not exist. 
% some Bertrand type curves do not appear. 
Moreover, the planar involutes and planar evolutes appear as the Bertrand type curves (Theorems \ref{tn-Bertrand-type} and \ref{nt-Bertrand-type}). 
For more details of properties of involutes and evolutes see \cite{Bruce-Giblin, Fukunaga-Takahashi-2015, Gibson, Gray, Nakatsuyama-Takahashi}.
\par
On the other hand, since the another curve may have singular points, we also consider curves with singular points. 
As smooth curves with singular points, it is useful to use the framed curves in the Euclidean space (cf. \cite{Honda-Takahashi-2016}). 
In \S 2, we also review the Bertrand and Mannheim curves of framed curves (cf. \cite{Honda-Takahashi-2020}). 
Then the same idea of Bertrand type curves is applied to the framed curves.
In \S 4,  we define and investigate Bertrand framed curves (Bertrand types of framed curves). 
We also clarify that the existence conditions of the Bertrand framed curves in all cases.
As a consequence, the involutes and circular evolutes of framed curves (cf. \cite{Honda-Takahashi-Preprint}) appear as the Bertrand framed curves (Theorems \ref{mu-nu1}, \ref{nu1-mu} and \ref{nu2-mu}). 
Therefore, it is useful to find new framed curves by using Bertrand framed curves. 
%(cf. \cite{Nakatsuyama-Takahashi}). 
\par
We shall assume throughout the whole paper that all maps and manifolds are $C^{\infty}$ unless the contrary is explicitly stated.

%%%%%%%%%%%%%%%%%%%%%%%%%%%%%%%%%%%%%%%%%%%%
\bigskip
\noindent
{\bf Acknowledgement}. 
%The authors would like to thank Professor Goo Ishikawa for his constant encouragement. 
The second author was supported by JSPS KAKENHI Grant Number 20K03573.

%%%%%%%%%%%%%%%%% Section 2 %%%%%%%%%%%%%%%%%
\section{Preliminaries}

We review the theories of Bertrand and Mannheim regular curves and framed curves.

Let $\R^3$ be the $3$-dimensional Euclidean space equipped with the inner product $\ba \cdot \bb = a_1 b_1 + a_2 b_2 + a_3 b_3$, 
where $\ba = (a_1, a_2, a_3)$ and $\bb = (b_1, b_2, b_3) \in \R^3$. 
The norm of $\ba$ is given by $\vert \ba \vert = \sqrt{\ba \cdot \ba}$ and the vector product is given by 
$$
\ba \times \bb={\rm det}
\left(
\begin{array}{ccc}
\be_1 & \be_2 & \be_3\\
a_1 & a_2 & a_3\\
b_1 & b_2 & b_3
\end{array}
\right)
$$
where $\{\be_1, \be_2, \be_3\}$ is the canonical basis of $\R^3$. 
Let $S^2$ be the unit sphere in $\R^3$, that is, $S^2=\{\ba \in \R^3| |\ba|=1\}$.
We denote the $3$-dimensional smooth manifold $\{(\ba,\bb) \in S^2 \times S^2| \ba \cdot \bb=0\}$ by $\Delta$. 

%%%%%
\subsection{Regular curves}

Let $I$ be an interval of $\R$ and let $\gamma:I \to \R^3$ be a regular space curve, that is, $\dot{\gamma}(t) \not=0$ for all $t \in I$, where $\dot{\gamma}(t)=(d\gamma/dt)(t)$. 
We say that $\gamma$ is {\it non-degenerate}, or $\gamma$ satisfies the {\it non-degenerate condition} if $\dot{\gamma}(t) \times \ddot{\gamma}(t) \not=0$ for all $t \in I$. 
%The non-degenerate condition is equivalent to the condition that the curvature of $\gamma$ is non-zero. 

If we take the arc-length parameter $s$, that is, $|\gamma'(s)|=1$ for all $s$, then the tangent vector, the principal normal vector and the bi-normal vector are given by
$$
\bt(s)=\gamma'(s), \ \bn(s)=\frac{\gamma''(s)}{|\gamma''(s)|}, \ \bb(s)=\bt(s) \times \bn(s),
$$
where $\gamma'(s)=(d\gamma/ds)(s)$. 
Then $\{\bt(s),\bn(s),\bb(s)\}$ is a moving frame of $\gamma(s)$ and we have the Frenet-Serret formula: 
$$
\left(
\begin{array}{c}
\bt'(s)\\
\bn'(s)\\
\bb'(s)
\end{array}
\right)
=
\left(
\begin{array}{ccc}
0&\kappa(s)&0\\
-\kappa(s)&0&\tau(s)\\
0&-\tau(s)&0
\end{array}
\right)
\left(
\begin{array}{c}
\bt(s)\\
\bn(s)\\
\bb(s)
\end{array}
\right),
$$
where 
$$
\kappa(s)=|\gamma''(s)|, \ \tau(s)=\frac{{\rm det}(\gamma'(s),\gamma''(s),\gamma'''(s))}{\kappa^2(s)}.
$$
If we take a general parameter $t$, then the tangent vector, the principal normal vector and the bi-normal vector are given by
$$
\bt(t)=\frac{\dot{\gamma}(t)}{|\dot{\gamma}(t)|}, \ \bn(t)=\bb(t) \times \bt(t), \ \bb(t)=\frac{\dot{\gamma}(t) \times \ddot{\gamma}(t)}{|\dot{\gamma}(t) \times \ddot{\gamma}(t)|}.
$$
Then $\{\bt(t),\bn(t),\bb(t)\}$ is a moving frame of $\gamma(t)$ and we have the Frenet-Serret formula: 
$$
\left(
\begin{array}{c}
\dot{\bt}(t)\\
\dot{\bn}(t)\\
\dot{\bb}(t)
\end{array}
\right)
=
\left(
\begin{array}{ccc}
0&|\dot{\gamma}(t)|\kappa(t)&0\\
-|\dot{\gamma}(t)|\kappa(t)&0&|\dot{\gamma}(t)|\tau(t)\\
0&-|\dot{\gamma}(t)|\tau(t)&0
\end{array}
\right)
\left(
\begin{array}{c}
\bt(t)\\
\bn(t)\\
\bb(t)
\end{array}
\right),
$$
where 
$$
\kappa(t)=\frac{|\dot{\gamma}(t) \times \ddot{\gamma}(t)|}{|\dot{\gamma}(t)|^3}, \ \tau(t)=\frac{{\rm det}(\dot{\gamma}(t),\ddot{\gamma}(t),\dddot{\gamma}(t))}{|\dot{\gamma}(t) \times \ddot{\gamma}(t)|^2}.
$$

Note that in order to define $\bt(t), \bn(t), \bb(t), \kappa(t)$ and $\tau(t)$, we assume that $\gamma$ is not only regular, but also non-degenerate.

%%%%%
\subsection{Bertrand and Mannheim non-degenerate curves}

Let $\gamma$ and $\overline{\gamma}:I \to \R^3$ be different non-degenerate curves. 

%%%%%
\begin{definition}\label{regular-Bertrand}{\rm 
We say that $\gamma$ and $\overline{\gamma}$ are {\it Bertrand mates} if there exists a smooth function $\lambda:I \to \R$ such that $\overline{\gamma}(t)=\gamma(t)+\lambda(t)\bn(t)$ and $\bn(t)=\pm \overline{\bn}(t)$ for all $t \in I$.
We also say that $\gamma:I \to \R^3$ is a {\it Bertrand curve} if there exists another non-degenerate curve $\overline{\gamma}:I \to \R^3$ such that $\gamma$ and $\overline{\gamma}$ are Bertrand mates.}
\end{definition}
%%%%%

In \cite{Honda-Takahashi-2020}, we investigated Bertrand curves of non-degenerate curves under the assumption $\tau(s) \not=0$ for all $s \in I$. 
However, we give an existence condition without this assumption.
By a parameter change, we may assume that $s$ is the arc-length parameter of $\gamma$.

%%%%%
\begin{theorem}\label{regular-Bertrand-equivalent}
Let $\gamma:I \to \R^3$ be non-degenerate with the arc-length parameter. 
\par
$(1)$ Suppose that there exists a point $s_0 \in I$ such that $\tau(s_0) \not=0$. 
Then $\gamma$ is a Bertrand curve if and only if there exists a non-zero constant $A$ and a constant $B$ such that  $A\kappa(s)+B\tau(s)=1$ and $\tau(s) (B \kappa(s)-A \tau(s)) \not=0$ for all $s \in I$. 
\par
$(2)$ Suppose that $\tau(s)=0$ for all $s \in I$. 
Then $\gamma$ is always a Bertrand curve. 
\end{theorem}
%%%%%
\demo
$(1)$ Suppose that $\gamma$ is a Bertrand curve. 
By differentiating $\overline{\gamma}(s)=\gamma(s)+\lambda(s)\bn(s)$, 
we have 
$$
|\dot{\overline{\gamma}}(s)|\overline{\bt}(s)=(1-\lambda(s)\kappa(s))\bt(s)+\lambda'(s)\bn(s)+\lambda(s)\tau(s)\bb(s).
$$ 
Since $\bn(s)=\pm \overline{\bn}(s)$, we have $\lambda'(s)=0$ for all $s \in I$. 
Therefore $\lambda$ is a constant. 
If $\lambda=0$, then $\overline{\gamma}(t)=\gamma(t)$ for all $t \in I$. 
Hence, $\lambda$ is a non-zero constant. 
We rewrite $\lambda$ as $A$. 
Note that $s$ is not the arc-length parameter of $\overline{\gamma}$. 
By differentiating $\overline{\gamma}(s)=\gamma(s)+A\bn(s)$, 
we have 
$$
|\dot{\overline{\gamma}}(s)|\overline{\bt}(s)=(1-A\kappa(s))\bt(s)+A\tau(s)\bb(s).
$$ 
If $\bn(s)=\overline{\bn}(s)$, there exists a smooth function $\theta:I \to \R$ such that
$$
\left(
\begin{array}{c}
\overline{\bb}(s)\\
\overline{\bt}(s)
\end{array}
\right)
=
\left(
\begin{array}{cc}
\cos \theta(s) &-\sin \theta(s)\\
\sin \theta(s) & \cos \theta(s)
\end{array}
\right)
\left(
\begin{array}{c}
{\bb}(s)\\
{\bt}(s)
\end{array}
\right).
$$
Then $|\dot{\overline{\gamma}}(s)|\sin \theta(s)=A \tau(s)$ and 
$|\dot{\overline{\gamma}}(s)|\cos \theta(s)=1-A \kappa(s)$. 
It follows that 
\begin{eqnarray}\label{equation1}
-A \tau(s) \cos \theta(s) +(1-A \kappa(s))\sin \theta(s)=0.
\end{eqnarray}
By differentiating $\overline{\bt}(s)=\sin \theta(s) \bb(s)+\cos \theta(s)\bt(s)$, we have
$$
|\dot{\overline{\gamma}}(s)|\overline{\kappa}(s)\overline{\bn}(s)=\theta'(s) \cos \theta(s) \bb(s)-\theta'(s) \sin \theta(s) \bt(s) +(-\tau(s) \sin \theta(s)+\kappa(s) \cos \theta(s))\bn(s). 
$$
Since $\bn(s)=\overline{\bn}(s)$, $\theta'(s)=0$ for all $s \in I$. 
Therefore $\theta$ is a constant. 
If $\tau(s)=0$ at some point $s \in I$, then $\sin \theta=0$ and hence $\tau(s)=0$ for all $s \in I$. 
It is a contradict of the condition $\tau(s_0) \not=0$. 
It follows that $\tau(s) \not=0$ for all $s \in I$ and $\sin \theta \not=0$.
By equation \eqref{equation1}, we have $A \kappa(s)+A (\cos\theta/\sin \theta)\tau(s)=1$. 
Hence, if we put $B=A \cos \theta/\sin \theta$, then  $A\kappa(s)+B\tau(s)=1$ for all $s \in I$. 
Moreover, 
$$
|\dot{\overline{\gamma}}(s)|\overline{\kappa}(s)=-\tau(s) \sin \theta +\kappa(s) \cos \theta =\frac{\sin \theta}{A}(-A\tau(s)+B\kappa(s)).
$$
Since $\overline{\kappa}(s) \not=0$, we have $B\kappa(s)-A\tau(s) \not=0$ for all $s \in I$.
\par
On the other hand, if $\bn(s)=-\overline{\bn}(s)$, there exists a smooth function $\theta:I \to \R$ such that
$$
\left(
\begin{array}{c}
\overline{\bb}(s)\\
\overline{\bt}(s)
\end{array}
\right)
=
\left(
\begin{array}{cc}
\cos \theta(s) &-\sin \theta(s)\\
\sin \theta(s) & \cos \theta(s)
\end{array}
\right)
\left(
\begin{array}{c}
{\bt}(s)\\
{\bb}(s)
\end{array}
\right).
$$
Then $|\dot{\overline{\gamma}}(s)|\sin \theta(s)=1-A \kappa(s)$ and 
$|\dot{\overline{\gamma}}(s)|\cos \theta(s)=A \tau(s)$. 
It follows that 
\begin{eqnarray}\label{equation2}
-A \tau(s) \sin \theta(s) +(1-A \kappa(s))\cos \theta(s)=0.
\end{eqnarray}
By differentiating $\overline{\bt}(s)=\sin \theta(s) \bt(s)+\cos \theta(s)\bb(s)$, we have
$$
|\dot{\overline{\gamma}}(s)|\overline{\kappa}(s)\overline{\bn}(s)=\theta'(s) \cos \theta(s) \bt(s)-\theta'(s) \sin \theta(s) \bb(s) +(\kappa(s) \sin \theta(s)-\tau(s) \cos \theta(s))\bn(s). 
$$
Since $\bn(s)=-\overline{\bn}(s)$, $\theta'(s)=0$ for all $s \in I$. 
Therefore $\theta$ is a constant. 
If $\tau(s)=0$ at some point $s \in I$, then $\cos \theta=0$ and hence $\tau(s)=0$ for all $s \in I$. 
It is a contradict of the condition $\tau(s_0) \not=0$. 
It follows that $\tau(s) \not=0$ for all $s \in I$ and $\cos \theta \not=0$.
By equation \eqref{equation2}, we have $A \kappa(s)+A (\sin \theta/\cos \theta)\tau(s)=1$. 
Hence, if we put $B=A \sin \theta/\cos \theta$, then  $A\kappa(s)+B\tau(s)=1$ for all $s \in I$. 
Moreover, 
$$
|\dot{\overline{\gamma}}(s)|\overline{\kappa}(s)=\kappa(s) \sin \theta+ \tau(s) \cos \theta 
=\frac{\cos \theta}{A}(-A\tau(s)+B\kappa(s)).
$$
Since $\overline{\kappa}(s) \not=0$, we have $B\kappa(s)-A\tau(s) \not=0$ for all $s \in I$. 
\par
Conversely, suppose there exists a non-zero constant $A$ and a constant $B$ such that  $A\kappa(s)+B\tau(s)=1$ and $\tau(s) (B \kappa(s)-A \tau(s)) \not=0$ for all $s \in I$. 
Set $\overline{\gamma}(s)=\gamma(s)+A\bn(s)$. 
By a direct calculation, we have
\begin{align*}
\dot{\overline{\gamma}}(s)&=|\dot{\overline{\gamma}}(s)|\overline{\bt}(s)=(1-A\kappa(s))\bt(s)+A\tau(s)\bb(s)=\tau(s)(B\bt(s)+A\bb(s)),\\
\ddot{\overline{\gamma}}(s)&=\tau'(s)(B\bt(s)+A\bb(s))+\tau(s)(B\kappa(s)-A\tau(s))\bn(s),\\
\dot{\overline{\gamma}}(s) \times \ddot{\overline{\gamma}}(s)&=\tau^2(s)B(B\kappa(s)-A\tau(s)) \bb(t)-\tau^2(s)A(B\kappa(s)-A\tau(s))\bt(s) \not=0.
\end{align*}
It follows that $\overline{\gamma}$ is a non-degenerate curve. 
Since $|\dot{\overline{\gamma}}(s)|=\sqrt{A^2+B^2} |\tau(s)|$, 
we have $\overline{\bt}(s)={\rm sgn}(\tau(s)) (1/\sqrt{A^2+B^2})(B \bt(s)+A\bb(s)),$ where ${\rm sgn}(\tau(s))=1$ if $\tau(s)>0$ and ${\rm sgn}(\tau(s))=-1$ if $\tau(s)<0$. 
By differentiating $\overline{\bt}(s)$, we have $|\dot{\overline{\gamma}}(s)|\overline{\kappa}(s)\overline{\bn}(s)={\rm sgn}(\tau(s)) (1/\sqrt{A^2+B^2})(B \kappa(s)-A \tau(s))\bn(s)$. 
By the assumption, we have $\bn(s)=\pm \overline{\bn}(s)$ for all $s \in I$. 
It follows that $\gamma$ and $\overline{\gamma}$ are Bertrand mates.
\par
$(2)$ Since $\kappa(s)>0$ for all $s \in I$, we can take a non-zero constant $\lambda$ such that $1-\lambda \kappa(s) \not= 0$ for all $s \in I$. 
In fact, if $\lambda$ is a negative constant, then the condition holds. 
Set $\overline{\gamma}(s)=\gamma(s)+\lambda \bn(s)$. 
By a direct calculation, we have 
\begin{align*}
\dot{\overline{\gamma}}(s)&=|\dot{\overline{\gamma}}(s)|\overline{\bt}(s)=(1-\lambda \kappa(s))\bt(s),\\
\ddot{\overline{\gamma}}(s)&=-\lambda \kappa'(s)\bt(s)+(1-\lambda \kappa(s))\kappa(s) \bn(s),\\
\dot{\overline{\gamma}}(s) \times \ddot{\overline{\gamma}}(s) &=(1-\lambda \kappa(s))^2\kappa(s)\bb(s) \not=0.
\end{align*}
It follows that $\overline{\gamma}$ is a non-degenerate curve. 
Since $|\dot{\overline{\gamma}}(s)|=|1-\lambda \kappa(s)|$ and $\overline{\bt}(s)=\pm \bt(s)$, we have $|\dot{\overline{\gamma}}(s)| \overline{\kappa}(s)\overline{\bn}(s)=\pm \kappa(s)\bn(s)$. 
By the assumption, we have $\bn(s)=\pm \overline{\bn}(s)$ for all $s \in I$. 
It follows that $\gamma$ and $\overline{\gamma}$ are Bertrand mates.
\enD
%%%%%
\begin{proposition}
Let $\gamma$ and $\overline{\gamma}:I\rightarrow\R^3$ be different non-degenerate curves. \par
$(1)$ Suppose that there exists a point $s_0 \in I$ such that $\tau(s_0) \not=0$ and $\gamma$ and $\overline{\gamma}$ are Bertrand mates with $\overline{\gamma}(s)=\gamma(s)+A\bn(s)$ and $A\kappa(s)+B\tau(s)=1$ for all $s\in I$, where $A$ is a non-zero constant and $B$ is a constant. Then the curvature $\overline{\kappa}$ and the torsion $\overline{\tau}$ of $\overline{\gamma}$ are given by
$$
\overline{\kappa}(s)=\frac{|B\kappa(s)-A\tau(s)|}{(A^2+B^2)|\tau(s)|}, \ \overline{\tau}(s)=\frac{1}{(A^2+B^2)\tau(s)}.
$$
\par
$(2)$ Suppose that $\tau(s)=0$ and $\gamma$ and $\overline{\gamma}$ are Bertrand mates with $\overline{\gamma}(s)=\gamma(s)+A\bn(s)$, where $A$ is non-zero constant and $1-A \kappa(s) \not=0$ for all $s\in I$. 
Then the curvature $\overline{\kappa}$ and the torsion $\overline{\tau}$ of $\overline{\gamma}$ are given by
$$
\overline{\kappa}(s)=\frac{\kappa(s)}{|1-A\kappa(s)|} , \ \overline{\tau}(s)=0.
$$
\end{proposition}
\demo
$(1)$ Since $\overline{\gamma}(s)=\gamma(s)+A\bn(s)$, we have 
\begin{align*}
\dot{\overline{\gamma}}(s)&=(1-A\kappa(s))\bt(s)+A\tau(s)\bb(s)=\tau(s)(B\bt(s)+A\bb(s)). 
\end{align*}
Therefore, we have
\begin{align*}
\ddot{\overline{\gamma}}(s) &=\tau'(s)(B\bt(s)+A\bb(s))+\tau(s)(B\kappa(s)-A\tau(s))\bn(s),\\
\dddot{\overline{\gamma}}(s) &=\tau''(s)(B\bt(s)+A\bb(s))+2\tau'(s)(B\kappa(s)-A\tau(s))\bn(s) \\
&\quad +\tau(s)(B\kappa'(s)-A\tau'(s))\bn(s)+\tau(s)(B\kappa(s)-A\tau(s))(-\kappa(s)\bt(s)+\tau(s)\bb(s)).
\end{align*}
Since
\begin{align*}
|\dot{\overline{\gamma}}(s)|&=|\tau(s)|\sqrt{A^2+B^2},\\
\dot{\overline{\gamma}}(s) \times \ddot{\overline{\gamma}}(s)&=\tau^2(s)(B\kappa(s)-A\tau(s))(B\bb(s)-A\bt(s)), \\
|\dot{\overline{\gamma}}(s) \times \ddot{\overline{\gamma}}(s)|&=\tau^2(s)|B\kappa(s)-A\tau(s)|\sqrt{A^2+B^2},\\
{\rm det}(\dot{\overline{\gamma}}(s), \ddot{\overline{\gamma}}(s), \dddot{\overline{\gamma}}(s))&=\tau^3(s)(B\kappa(s)-A\tau(s))^2,
\end{align*}
we have the curvature $\overline{\kappa}$ and the torsion $\overline{\tau}$ of $\overline{\gamma}$ as
\begin{align*}
\overline{\kappa}(s)&=\frac{|\dot{\overline{\gamma}}(s) \times \ddot{\overline{\gamma}}(s)|}{|\dot{\overline{\gamma}}(s)|^3}=\frac{|B\kappa(s)-A\tau(s)|}{(A^2+B^2)|\tau(s)|}, \\ 
\overline{\tau}(s)&=\frac{{\rm det}(\dot{\overline{\gamma}}(s), \ddot{\overline{\gamma}}(s), \dddot{\overline{\gamma}}(s))}{|\dot{\overline{\gamma}}(s) \times \ddot{\overline{\gamma}}(s)|^2}=\frac{1}{(A^2+B^2)\tau(s)}.
\end{align*}
\par
$(2)$ Since $\overline{\gamma}(s)=\gamma(s)+A\bn(s)$, we have 
\begin{align*}
\dot{\overline{\gamma}}(s)&=(1-A\kappa(s))\bt(s), \\
\ddot{\overline{\gamma}}(s)&=-A\kappa'(s)\bt(s)+(1-A\kappa(s))\kappa(s)\bn(s),\\
\dddot{\overline{\gamma}}(s)&=\left(-A\kappa''(s)-(1-A\kappa(s))\kappa^2(s)\right)\bt(s)+\left(\kappa'(s)-3A\kappa(s)\kappa'(s)\right)\bn(t).
\end{align*}
Since
\begin{align*}
|\dot{\overline{\gamma}}(s)|&=|1-A\kappa(s)|,\\
\dot{\overline{\gamma}}(s) \times \ddot{\overline{\gamma}}(s)&=\kappa(s)(1-A\kappa(s))^2\bb(s),\\
|\dot{\overline{\gamma}}(s) \times \ddot{\overline{\gamma}}(s)|&=\kappa(s)(1-A\kappa(s))^2,\\
{\rm det}(\dot{\overline{\gamma}}(s), \ddot{\overline{\gamma}}(s), \dddot{\overline{\gamma}}(s))&=0,
\end{align*}
we have the curvature $\overline{\kappa}$ and the torsion $\overline{\tau}$ of $\overline{\gamma}$ as
\begin{align*}
\overline{\kappa}(s)=\frac{|\dot{\overline{\gamma}}(s) \times \ddot{\overline{\gamma}}(s)|}{|\dot{\overline{\gamma}}(s)|^3}=\frac{\kappa(s)}{|1-A\kappa(s)|}, \ 
\overline{\tau}(s)=\frac{{\rm det}(\dot{\overline{\gamma}}(s), \ddot{\overline{\gamma}}(s), \dddot{\overline{\gamma}}(s))}{|\dot{\overline{\gamma}}(s) \times \ddot{\overline{\gamma}}(s)|^2}=0.
\end{align*}
\enD

%%%%%
\begin{definition}\label{regular-Mannheim}{\rm 
We say that $\gamma$ and $\overline{\gamma}$ are {\it Mannheim mates} if there exists a smooth function $\lambda:I \to \R$ such that $\overline{\gamma}(t)=\gamma(t)+\lambda(t)\bn(t)$ and $\bn(t)=\pm \overline{\bb}(t)$ for all $t \in I$.
We also say that $\gamma:I \to \R^3$ is a {\it Mannheim curve} if there exists another non-degenerate curve $\overline{\gamma}:I \to \R^3$ such that $\gamma$ and $\overline{\gamma}$ are Mannheim mates.}
\end{definition}
%%%%%

In \cite{Honda-Takahashi-2020}, we also investigated Mannheim curves of non-degenerate curves under the assumption $\tau(s) \not=0$ for all $s \in I$. 
However, the assumption does not need. 
We give a result and proof in detail.

%%%%%
\begin{theorem}\label{regular-Mannheim-equivalent}
Let $\gamma:I \to \R^3$ be non-degenerate with the arc-length parameter. 
Then $\gamma$ is a Mannheim curves if and only if there exists a non-zero constant $A$ such that  $A(\kappa^2(s)+\tau^2(s))=\kappa(s)$ and $\tau(s)(\kappa(s)\tau'(s)-\kappa'(s)\tau(s)) \not=0$ for all $s \in I$.
\end{theorem}
%%%%%
\demo 
Suppose that $\gamma$ is a Mannheim curve. 
There exists a smooth function $\lambda:I \to \R$ such that $\overline{\gamma}(t)=\gamma(t)+\lambda(t)\bn(t)$ and $\bn(t)=\pm \overline{\bb}(t)$ for all $t \in I$. 
By differentiating $\overline{\gamma}(s)=\gamma(s)+\lambda(s)\bn(s)$, 
we have 
$
|\dot{\overline{\gamma}}(s)|\overline{\bt}(s)=(1-\lambda(s)\kappa(s))\bt(s)+\lambda'(s)\bn(s)+\lambda(s)\tau(s)\bb(s).
$ 
Since $\bn(s)=\pm \overline{\bb}(s)$, we have $\lambda'(s)=0$ for all $s \in I$. Therefore $\lambda$ is a constant. If $\lambda=0$, then $\overline{\gamma}(t)=\gamma(t)$ for all $t \in I$. 
Hence, $\lambda$ is a non-zero constant.
We rewrite $\lambda$ as $A$.
Note that $s$ is not the arc-length parameter of $\overline{\gamma}$. 
By differentiating $\overline{\gamma}(s)=\gamma(s)+A\bn(s)$, 
we have 
$$
|\dot{\overline{\gamma}}(s)|\overline{\bt}(s)=(1-A\kappa(s))\bt(s)+A\tau(s)\bb(s).
$$ 
If $\bn(s)=\overline{\bb}(s)$, there exists a smooth function $\theta:I \to \R$ such that
$$
\left(
\begin{array}{c}
\overline{\bt}(s)\\
\overline{\bn}(s)
\end{array}
\right)
=
\left(
\begin{array}{cc}
\cos \theta(s) &-\sin \theta(s)\\
\sin \theta(s) & \cos \theta(s)
\end{array}
\right)
\left(
\begin{array}{c}
{\bb}(s)\\
{\bt}(s)
\end{array}
\right).
$$
Then $|\dot{\overline{\gamma}}(s)|\cos \theta(s)=A \tau(s)$ and 
$-|\dot{\overline{\gamma}}(s)|\sin \theta(s)=1-A \kappa(s)$. 
It follows that 
\begin{eqnarray}\label{equation3}
A \tau(s) \sin \theta(s)+(1-A \kappa(s))\cos \theta(s)=0.
\end{eqnarray}
By differentiating $\overline{\bt}(s)=\cos \theta(s) \bb(s)-\sin \theta(s)\bt(s)$, we have
$$
|\dot{\overline{\gamma}}(s)|\overline{\kappa}(s)\overline{\bn}(s)=-\theta'(s) \sin \theta(s) \bb(s)-\theta'(s) \cos \theta(s) \bt(s) -(\tau(s) \cos \theta(s)+\kappa(s) \sin \theta(s))\bn(s). 
$$
Since $\bn(s)=\overline{\bb}(s)$, 
\begin{eqnarray}\label{equation4}
\tau(s) \cos \theta(s)+\kappa(s) \sin \theta(s)=0
\end{eqnarray} for all $s \in I$.  
If $\tau(s)=0$ at a point $s \in I$, then $\cos \theta(s)=0$ and $\sin \theta(s)=0$ by \eqref{equation4}.
Hence $\tau(s) \not=0$ for all $s \in I$. 
It follows that $\cos \theta(s) \not=0$ and $\sin \theta(s) \not=0$.
Since $\overline{\bn}(s)=\sin \theta(s) \bb(s)+\cos \theta(s) \bt(s)$, we have 
$|\dot{\overline{\gamma}}(s)| \overline{\kappa}(s)=-\theta'(s)$. 
By equations \eqref{equation3} and \eqref{equation4}, we have
\begin{eqnarray*}
A(\kappa^2(s)+\tau^2(s))=\kappa(s).
\end{eqnarray*}
By differentiating \eqref{equation4}, we have 
$$
-\tau(s) \theta'(s) \sin \theta(s) +\tau'(s) \cos \theta(s)+\kappa(s) \theta'(s)\cos \theta(s)+\kappa'(s) \sin \theta(s)=0.
$$
Hence $\theta'(s)=(-\kappa(s)\tau'(s)+\kappa'(s)\tau(s))/(\kappa^2(s)+\tau^2(s))$. 
Since $|\dot{\overline{\gamma}}(s)|\overline{\kappa}(s) >0$, we have 
$\kappa(s) \tau'(s)-\kappa'(s) \tau(s) >0$ for all $s \in I$. 
\par
On the other hand, if $\bn(s)=-\overline{\bb}(s)$, there exists a smooth function $\theta:I \to \R$ such that
$$
\left(
\begin{array}{c}
\overline{\bt}(s)\\
\overline{\bn}(s)
\end{array}
\right)
=
\left(
\begin{array}{cc}
\cos \theta(s) &-\sin \theta(s)\\
\sin \theta(s) & \cos \theta(s)
\end{array}
\right)
\left(
\begin{array}{c}
{\bt}(s)\\
{\bb}(s)
\end{array}
\right).
$$
Then $|\dot{\overline{\gamma}}(s)|\cos \theta(s)=1-A \kappa(s)$ and 
$-|\dot{\overline{\gamma}}(s)|\sin \theta(s)=A \tau(s)$. 
It follows that 
\begin{eqnarray}\label{equation5}
A \tau(s) \sin \theta(s)+(1-A \kappa(s))\cos \theta(s)=0.
\end{eqnarray}
By differentiating $\overline{\bt}(s)=\cos \theta(s) \bt(s)-\sin \theta(s)\bb(s)$, we have
$$
|\dot{\overline{\gamma}}(s)|\overline{\kappa}(s)\overline{\bn}(s)=-\theta'(s) \sin \theta(s) \bt(s)-\theta'(s) \cos \theta(s) \bb(s) -(\tau(s) \sin \theta(s) +\kappa(s) \cos \theta(s))\bn(s). 
$$
Since $\bn(s)=-\overline{\bb}(s)$, 
\begin{eqnarray}\label{equation6}
\tau(s) \sin \theta(s) + \kappa(s) \cos \theta(s)=0
\end{eqnarray} for all $s \in I$.  
If $\tau(s)=0$ at a point $s \in I$, then $\sin \theta(s)=0$ and $\cos \theta(s)=0$ by \eqref{equation6}.
Hence $\tau(s) \not=0$ for all $s \in I$. 
It follows that $\sin \theta(s) \not=0$ and $\cos \theta(s) \not=0$.
Since $\overline{\bn}(s)=\sin \theta(s) \bt(s)+\cos \theta(s) \bb(s)$, we have 
$|\dot{\overline{\gamma}}(s)| \overline{\kappa}(s)=-\theta'(s)$. 
By equations \eqref{equation5} and \eqref{equation6}, we have
\begin{eqnarray*}
A(\kappa^2(s)+\tau^2(s))=\kappa(s).
\end{eqnarray*}
By differentiating \eqref{equation6}, we have 
$$
\tau(s) \theta'(s) \cos \theta(s) +\tau'(s) \sin \theta(s)-\kappa(s) \theta'(s)\sin \theta(s)+\kappa'(s) \cos \theta(s)=0.
$$
Hence $\theta'(s)=(\kappa(s)\tau'(s)-\kappa'(s)\tau(s))/(\kappa^2(s)+\tau^2(s))$. 
Since $|\dot{\overline{\gamma}}(s)|\overline{\kappa}(s) >0$, we have 
$\kappa(s) \tau'(s)-\kappa'(s) \tau(s) <0$ for all $s \in I$. 
\par
Conversely, suppose that there exists a non-zero constant $A$ such that  $A(\kappa^2(s)+\tau^2(s))=\kappa(s)$ and $\tau(s)(\kappa(s)\tau'(s)-\kappa'(s)\tau(s)) \not=0$ for all $s \in I$. 
Set $\overline{\gamma}(s)=\gamma(s)+A\bn(s)$. 
By a direct calculation, we have
\begin{align*}
\dot{\overline{\gamma}}(s)&=|\dot{\overline{\gamma}}(s)| \overline{\bt}(t)=(1-A\kappa(s))\bt(s)+A\tau(s)\bb(s)=A\frac{\tau(s)}{\kappa(s)}(\tau(s)\bt(s)+\kappa(s)\bb(s)),\\
\ddot{\overline{\gamma}}(s)&=A \left(\frac{\tau(s)}{\kappa(s)}\right)'(\tau(s)\bt(s)+\kappa(s)\bb(s))+A\frac{\tau(s)}{\kappa(s)}(\tau'(s)\bt(s)+\kappa'(s)\bb(s)),\\
\dot{\overline{\gamma}}(s) \times \ddot{\overline{\gamma}}(s) &=A^2 \frac{\tau^2(s)}{\kappa^2(s)}(\kappa(s)\tau'(s)-\kappa'(s)\tau(s))\bn(s) \not=0.
\end{align*}
It follows that $\overline{\gamma}$ is a non-degenerate curve. 
Since $\ddot{\overline{\gamma}}(s)=({d}/{ds})(|\dot{\overline{\gamma}}(s)|)\overline{\bt}(s)+|\dot{\overline{\gamma}}(s)|^2 \overline{\kappa}(s) \overline{\bn}(s)$, we have $|\dot{\overline{\gamma}}(s)|^3 \overline{\kappa}(s) \overline{\bb}(s)=A^2(\tau^2 (s)/\kappa^2(s))(\kappa(s)\tau'(s)-\kappa'(s)\tau(s))\bn(s)$. 
By the assumption, we have $\bn(s)=\pm \overline{\bb}(s)$. 
It follows that $\gamma$ and $\overline{\gamma}$ are Mannheim mates.
\enD
%%%%%
\begin{proposition}
Let $\gamma$ and $\overline{\gamma}:I\rightarrow\R^3$ be different non-degenerate curves. Under the same assumptions in Theorem \ref{regular-Mannheim-equivalent}, suppose that $\gamma$ and $\overline{\gamma}$ are Mannheim mates with $\overline{\gamma}(s)=\gamma(s)+A\bn(s)$ for all $s\in I$, where $A$ is a non-zero constant. Then the curvature $\overline{\kappa}$ and the torsion $\overline{\tau}$ of $\overline{\gamma}$ are given by
$$
\overline{\kappa}(s)=\frac{\kappa(s)|\kappa(s)\tau'(s)-\kappa'(s)\tau(s)|}{|A\tau(s)|(\tau^2(s)+\kappa^2(s))^{\frac{3}{2}}} , \ \overline{\tau}(s)=\frac{\kappa^2(s)+\tau^2(s)}{\tau(s)}.
$$
\end{proposition}
\demo
Since $\overline{\gamma}(s)=\gamma(s)+A\bn(s)$, we have 
\begin{align*}
\dot{\overline{\gamma}}(s)&=(1-A\kappa(s))\bt(s)+A\tau(s)\bb(s)=A\frac{\tau(s)}{\kappa(s)}(\tau(s)\bt(s)+\kappa(s)\bb(s)). 
\end{align*}
Therefore, we have
\begin{align*}
\ddot{\overline{\gamma}}(s) &=A\left(\frac{\tau(s)}{\kappa(s)}\right)'(\tau(s)\bt(s)+\kappa(s)\bb(s))+A\frac{\tau(s)}{\kappa(s)}(\tau'(s)\bt(s)+\kappa'(s)\bb(s)),\\
\dddot{\overline{\gamma}}(s) &=A\left(\frac{\tau(s)}{\kappa(s)}\right)''(\tau(s)\bt(s)+\kappa(s)\bb(s))+2A\left(\frac{\tau(s)}{\kappa(s)}\right)'(\tau'(s)\bt(s)+\kappa'(s)\bb(s)),\\
&\quad +A\frac{\tau(s)}{\kappa(s)}(\tau''(s)\bt(s)+\kappa(s)\tau'(s)\bn(s)+\kappa''(s)\bb(s)-\kappa'(s)\tau(s)\bn(s)).
\end{align*}
Since
\begin{align*}
|\dot{\overline{\gamma}}(s)| &= \frac{|A\tau(s)|}{\kappa(s)}\sqrt{\tau^2(s)+\kappa^2(s)},\\
\dot{\overline{\gamma}}(s) \times \ddot{\overline{\gamma}}(s) &= A^2\left(\frac{\tau(s)}{\kappa(s)}\right)^2(\kappa(s)\tau'(s)-\kappa'(s)\tau(s))\bn(s),\\
|\dot{\overline{\gamma}}(s) \times \ddot{\overline{\gamma}}(s)| &= A^2\left(\frac{\tau(s)}{\kappa(s)}\right)^2|\kappa(s)\tau'(s)-\kappa'(s)\tau(s)|,\\
{\rm det}(\dot{\overline{\gamma}}(s), \ddot{\overline{\gamma}}(s), \dddot{\overline{\gamma}}(s)) &= A^3\left(\frac{\tau(s)}{\kappa(s)}\right)^3(\kappa(s)\tau'(s)-\kappa'(s)\tau(s))^2,
\end{align*}
we have the curvature and the torsion as
\begin{align*}
\overline{\kappa}(s)&=\frac{|\dot{\overline{\gamma}}(s) \times \ddot{\overline{\gamma}}(s)|}{|\dot{\overline{\gamma}}(s)|^3}=\frac{\kappa(s)|\kappa(s)\tau'(s)-\kappa'(s)\tau(s)|}{|A\tau(s)|(\tau^2(s)+\kappa^2(s))^{\frac{3}{2}}}, \\ 
\overline{\tau}(s)&=\frac{{\rm det}(\dot{\overline{\gamma}}(s), \ddot{\overline{\gamma}}(s), \dddot{\overline{\gamma}}(s))}{|\dot{\overline{\gamma}}(s) \times \ddot{\overline{\gamma}}(s)|^2}=\frac{\kappa^2(s)+\tau^2(s)}{\tau(s)}.
\end{align*}
\enD
%%%%%

\subsection{Framed curves}

A framed curve in the $3$-dimensional Euclidean space is a smooth space curve with a moving frame, in detail see \cite{Honda-Takahashi-2016}. 
%%%%%
\begin{definition}\label{framed.curve}{\rm
We say that $(\gamma,\nu_1,\nu_2):I \rightarrow \mathbb{R}^3 \times \Delta$ is a {\it framed curve} if $\dot{\gamma}(t) \cdot \nu_1(t)=0$ and $\dot{\gamma}(t) \cdot \nu_2(t)=0$ for all $t \in I$. 
We say that $\gamma:I \to \R^3$ is a {\it framed base curve} if there exists $(\nu_1,\nu_2):I \to \Delta$ such that $(\gamma,\nu_1,\nu_2)$ is a framed curve. 
}
\end{definition}
%%%%%

We denote $\mu(t) = \nu_1(t) \times \nu_2(t)$. 
Then $\{ \nu_1(t),\nu_2(t),\bmu(t) \}$ is a moving frame along the framed base curve $\gamma(t)$ in $\R^3$ and we have the Frenet type formula,
$$
\left(
\begin{array}{c}
\dot{\nu_1}(t)\\
\dot{\nu_2}(t)\\
\dot{\mu}(t)
\end{array} \right)=
\left(
\begin{array}{ccc}
0 & \ell(t) & m(t)\\
-\ell(t) & 0 & n(t)\\
-m(t) & -n(t) & 0
\end{array}\right)
\left(
\begin{array}{c}
\nu_1(t)\\
\nu_2(t)\\
\mu(t)
\end{array}\right), 
\ \dot{\gamma}(t)=\alpha(t)\mu(t),
$$
where $\ell(t) = \dot{\nu_1}(t) \cdot \nu_2(t)$, $m(t) = \dot{\nu_1}(t) \cdot \mu(t), n(t) = \dot{\nu_2}(t) \cdot \mu(t)$ and $\alpha(t)=\dot{\gamma}(t) \cdot \mu(t)$. 
We call the mapping $(\ell,m,n,\alpha)$ {\it the curvature of the framed curve} $(\gamma,\nu_1,\nu_2)$. 
Note that $t_0$ is a singular point of $\gamma$ if and only if $\alpha(t_0) = 0$. 
%%%%%
\begin{definition}
{\rm 
Let $(\gamma,\nu_1,\nu_2)$ and $(\widetilde{\gamma},\widetilde{\nu}_1,\widetilde{\nu}_2):I \rightarrow \mathbb{R}^3 \times \Delta$ be framed curves. 
We say that $(\gamma,\nu_1,\nu_2)$ and $(\widetilde{\gamma},\widetilde{\nu}_1,\widetilde{\nu}_2)$ are {\it congruent as framed curves} if there exist a constant rotation $A \in SO(3)$ and a translation $\ba \in \mathbb{R}^3$ such that $\widetilde{\gamma}(t) = A(\gamma(t)) +\ba$, $\widetilde{\nu_1}(t) = A(\nu_1(t))$ and $\widetilde{\nu_2}(t) = A(\nu_2(t))$ for all $t \in I$.
}
\end{definition}
%%%%%
\par 
We gave the existence and uniqueness theorems for framed curves in terms of the curvatures in \cite{Honda-Takahashi-2016}, also see \cite{Fukunaga-Takahashi-2017}.  
%%%%%
\begin{theorem}[Existence Theorem for framed curves]\label{existence.framed}
Let $(\ell,m,n,\alpha):I \rightarrow \mathbb{R}^4$ be a smooth mapping. 
Then, there exists a framed curve $(\gamma,\nu_1,\nu_2):I \to \R^3
 \times \Delta$ whose curvature is given by $(\ell,m,n,\alpha)$. 
\end{theorem}
%%%%%
\begin{theorem}[Uniqueness Theorem for framed curves]\label{uniqueness.framed}
Let $(\gamma,\nu_1,\nu_2)$ and $(\widetilde{\gamma},\widetilde{\nu}_1,\widetilde{\nu}_2):I \to \R^3 \times \Delta$ be framed curves with curvatures $(\ell,m,n,\alpha)$ and $(\widetilde{\ell},\widetilde{m},\widetilde{n},\widetilde{\alpha})$, respectively. 
Then $(\gamma,\nu_1,\nu_2)$ and $(\widetilde{\gamma},\widetilde{\nu}_1,\widetilde{\nu}_2)$ are congruent as framed curves if and only if the curvatures $(\ell, m, n, \alpha)$ and $(\widetilde{\ell}, \widetilde{m}, \widetilde{n}, \widetilde{\alpha})$ coincide. 
\end{theorem}
%%%%%

Let $(\gamma,\nu_1,\nu_2): I \to \R^3 \times \Delta$ be a framed curve with the curvature $(\ell,m,n,\alpha)$. 
For the normal plane of $\gamma(t)$, spanned by $\nu_1(t)$ and $\nu_2(t)$, there is some ambient of framed curves similarly to the case of the Bishop frame of a regular space curve (cf. \cite{Bishop}). 
We define $(\widetilde{\nu}_1(t), \widetilde{\nu}_2(t)) \in \Delta_2$ by 
$$
\left(\begin{array}{c}
\widetilde{\nu}_1(t)\\
\widetilde{\nu}_2(t)
\end{array}\right)
=
\left(\begin{array}{cc}
\cos \theta(t) & -\sin \theta(t)\\
\sin \theta(t) & \cos \theta(t)
\end{array}\right)
\left(\begin{array}{c}
\nu_1(t)\\
\nu_2(t)
\end{array}\right),
$$
where $\theta(t)$ is a smooth function.
Then $(\gamma,\widetilde{\nu}_1,\widetilde{\nu}_2): I \to \R^3 \times \Delta$ is also a framed curve and $\widetilde{\mu}(t)=\mu(t).$
By a direct calculation, we have
\begin{align*}
\dot{\widetilde{\nu}}_1(t) &= (\ell(t)-\dot{\theta}(t))\sin \theta(t) \nu_1(t)+(\ell(t)-\dot{\theta}(t))\cos \theta(t) \nu_2(t) \\
&\quad +(m(t)\cos \theta(t)-n(t)\sin \theta(t))\mu(t), \\
\dot{\widetilde{\nu}}_2(t) &= -(\ell(t)-\dot{\theta}(t))\cos \theta(t) \nu_1(t)+(\ell(t)-\dot{\theta}(t))\sin \theta(t) \nu_2(t) \\
&\quad +(m(t)\sin \theta(t)+n(t)\cos \theta(t))\mu(t). 
\end{align*}
If we take a smooth function $\theta:I \rightarrow \mathbb{R}$ which satisfies $\dot{\theta}(t)=\ell(t)$, then we call the frame $\{\widetilde{\nu}_1(t),\widetilde{\nu}_2(t),\mu(t)\}$ an {\it adapted frame} along $\gamma(t)$.
It follows that the Frenet-Serret type formula is given by 
\begin{eqnarray}\label{adapted-frame}
\left(
\begin{array}{c}
\dot{\widetilde{\nu}}_1(t)\\
\dot{\widetilde{\nu}}_2(t)\\
\dot{\mu}(t)
\end{array} \right)=
\left(
\begin{array}{ccc}
0 & 0 & \widetilde{m}(t)\\
0 & 0 & \widetilde{n}(t)\\
-\widetilde{m}(t) & -\widetilde{n}(t) & 0
\end{array}\right)
\left(
\begin{array}{c}
\widetilde{\nu}_1(t)\\
\widetilde{\nu}_2(t)\\
\mu(t)
\end{array}\right),
\end{eqnarray}
where $\widetilde{m}(t)$ and $\widetilde{n}(t)$ are given by
\begin{eqnarray*}\label{adapted.curvature}
\left(\begin{array}{c}
\widetilde{m}(t)\\
\widetilde{n}(t)
\end{array}\right)
=
\left(\begin{array}{cc}
\cos \theta(t) & -\sin \theta(t)\\
\sin \theta(t) & \cos \theta(t)
\end{array}\right)
\left(\begin{array}{c}
m(t)\\
n(t)
\end{array}\right).
\end{eqnarray*}
%%%%%
\par
By a direct calculation, we have the following (cf. \cite{Honda-Takahashi-2016}).
%%%%%
\begin{proposition}\label{change-frame}
Suppose that $(\gamma,\nu_1,\nu_2):I \to \R^3 \times \Delta$ is a framed curve with curvature $(\ell,m,n,\alpha)$.
Then $(\gamma,\nu_2,\nu_1):I \to \R^3 \times \Delta$ is also a framed curve with curvature $(-\ell,-n,-m,-\alpha)$.
\end{proposition}
%%%%%
\subsection{Bertrand and Mannheim curves of framed curves}

Let $(\gamma,\nu_1,\nu_2)$ and $(\overline{\gamma},\overline{\nu}_1,\overline{\nu}_2):I \to \R^3 \times \Delta$ be framed curves with the curvature $(\ell,m,n,\alpha)$ and $(\overline{\ell},\overline{m},\overline{n},\overline{\alpha})$, respectively. 
Suppose that $\gamma$ and $\overline{\gamma}$ are different curves. 
In \cite{Honda-Takahashi-2020}, we define Bertrand and Mannheim curves of framed curves, and give a characterization of the Bertrand and Mannheim curves.

%%%%%
\begin{definition}\label{Bertrand.framed}{\rm 
We say that framed curves $(\gamma,\nu_1,\nu_2)$ and $(\overline{\gamma},\overline{\nu}_1,\overline{\nu}_2)$ are {\it Bertrand mates} (or, $(\nu_1,\overline{\nu}_1)$-mates) if there exists a smooth function $\lambda:I \to \R$ such that $\overline{\gamma}(t)=\gamma(t)+\lambda(t)\nu_1(t)$ and $\nu_1(t)=\overline{\nu}_1(t)$ for all $t \in I$. 
We also say that $(\gamma,\nu_1,\nu_2):I \to \R^3 \times \Delta$ is a {\it Bertrand curve} if there exists a framed curve  $(\overline{\gamma},\overline{\nu}_1,\overline{\nu}_2):I \to \R^3 \times \Delta$ such that $(\gamma,\nu_1,\nu_2)$ and $(\overline{\gamma},\overline{\nu}_1,\overline{\nu}_2)$ are Bertrand mates.
}
\end{definition}
%%%%%

%%%%%
\begin{theorem}\label{framed-Bertrand-equivalent}
Let $(\gamma,\nu_1,\nu_2):I \to \R^3 \times \Delta$ be a framed curve with the curvature $(\ell,m,n,\alpha)$.
Then $(\gamma,\nu_1,\nu_2)$ is a Bertrand curve if and only if there exists a non-zero constant $\lambda$ and a smooth function $\theta:I \to \R$ such that 
$\lambda \ell(t) \cos \theta(t)-(\alpha(t)+\lambda m(t))\sin \theta(t)=0
$ for all $t \in I$.
\end{theorem}
%%%%%
\begin{proposition}\label{framed-Bertrand_curvature}
Suppose that $(\gamma,\nu_1,\nu_2)$ and $(\overline{\gamma},\overline{\nu}_1,\overline{\nu}_2): I \to \R^3 \times \Delta$ are {Bertrand mates}, where $\overline{\gamma}(t)=\gamma(t)+\lambda\nu_1(t), \overline{\nu}_1(t)=\sin\theta(t)\nu_2(t)+\cos\theta(t)\mu(t), \overline{\nu}_2(t)=\nu_1(t)$ and $\theta:I \to \R$ is a smooth function. Then the curvature $(\overline{\ell}, \overline{m}, \overline{n}, \overline{\alpha})$ of $(\overline{\gamma},\overline{\nu}_1,\overline{\nu}_2)$ is given by
\begin{align*}
  \overline{\ell}(t)&=\ell(t)\cos\theta(t)-m(t)\sin\theta(t), \\
  \overline{m}(t)&=\ell(t)\sin\theta(t)+m(t)\cos\theta(t), \\
  \overline{n}(t)&=n(t)-\dot{\theta}(t), \\
  \overline{\alpha}(t)&=\lambda\ell(t)\sin\theta(t)+(\alpha(t)+\lambda m(t))\cos\theta(t).
\end{align*}
\end{proposition}

%%%%%
\begin{definition}\label{Mannheim.framed}{\rm 
We say that framed curves $(\gamma,\nu_1,\nu_2)$ and $(\overline{\gamma},\overline{\nu}_1,\overline{\nu}_2)$ are {\it Mannheim mates} (or, $(\nu_1,\overline{\nu}_2)$-mates) if there exists a smooth function $\lambda:I \to \R$ such that $\overline{\gamma}(t)=\gamma(t)+\lambda(t)\nu_1(t)$ and $\nu_1(t)=\overline{\nu}_2(t)$ for all $t \in I$. 
We also say that $(\gamma,\nu_1,\nu_2):I \to \R^3 \times \Delta$ is a {\it Mannheim curve} if there exists a framed curve  $(\overline{\gamma},\overline{\nu}_1,\overline{\nu}_2):I \to \R^3 \times \Delta$ such that $(\gamma,\nu_1,\nu_2)$ and $(\overline{\gamma},\overline{\nu}_1,\overline{\nu}_2)$ are Mannheim mates.
}
\end{definition}
%%%%%

%%%%%
\begin{theorem}\label{framed-Mannheim-equivalent}
Let $(\gamma,\nu_1,\nu_2):I \to \R^3 \times \Delta$ be a framed curve with the curvature $(\ell,m,n,\alpha)$.
Then $(\gamma,\nu_1,\nu_2)$ is a Mannheim curve if and only if there exists a non-zero constant $\lambda$ and a smooth function $\theta:I \to \R$ such that 
$\lambda \ell(t) \sin \theta(t)+(\alpha(t)+\lambda m(t))\cos \theta(t)=0
$ for all $t \in I$.
\end{theorem}
%%%%%

As a difference between non-degenerate regular space curves and framed curves, 
we have a relation between Bertrand and Mannheim curves of framed curves.
%%%%%
\begin{theorem}\label{relation}
Let $(\gamma,\nu_1,\nu_2):I \to \R^3 \times \Delta$ be a framed curve with the curvature $(\ell,m,n,\alpha)$. 
Then $(\gamma,\nu_1,\nu_2)$ is a Bertrand curve if and only if $(\gamma,\nu_1,\nu_2)$ is a Mannheim curve.
\end{theorem}
%%%%%
\begin{remark}\label{relation-remark}{\rm
By definition, if $(\gamma,\nu_1,\nu_2)$ and $(\overline{\gamma},\overline{\nu}_1,\overline{\nu}_2)$ are Bertrand mates, then $\overline{\gamma}=\gamma+\lambda \nu_1$ and $\nu_1=\overline{\nu}_1$.
Since $(\overline{\gamma},\overline{\nu}_2,\overline{\nu}_1)$ is also a framed curve, $(\gamma,\nu_1,\nu_2)$ and $(\overline{\gamma},\overline{\nu}_2,\overline{\nu}_1)$ are Mannheim mates and vice versa. 
Hence, $(\gamma,\nu_1,\nu_2)$ is a Bertrand curve if and only if $(\gamma,\nu_1,\nu_2)$ is a Mannheim curve.
}
\end{remark}
%%%%%

%%%%%%%%%%%%%%%%% Section 3 %%%%%%%%%%%%%%%%%
\section{Bertrand types of non-degenerate curves}

Let $\gamma$ and $\overline{\gamma}:I \to \R^3$ be non-degenerate curves. 
We consider a space curve whose tangent (or, principal normal, bi-normal) line is the same as the tangent (or, principal normal, bi-normal) line of another curve, respectively. 

%%%%%
\begin{definition}\label{Bertrand-type-regular}{\rm
We say that $\gamma$ and $\overline{\gamma}$ are {\it $(\bv,\overline{\bw})$-mates} if there exists a smooth function $\lambda:I \to \R$ with $\lambda \not\equiv 0$ such that $\overline{\gamma}(t)=\gamma(t)+\lambda(t)\bv(t)$ and $\bv(t)=\pm \overline{\bw}(t)$ for all $t \in I$, 
where $\bv$ and $\bw$ are $\bt, \bn$ or $\bb$. 
We also say that $\gamma$ is a {\it $(\bv,\overline{\bw})$-Bertrand type curve} (or, {\it $(\bv,\overline{\bw})$-Bertrand-Mannheim type curve}) if there exists another non-degenerate regular curve $\overline{\gamma}$ such that $\gamma$ and $\overline{\gamma}$ are $(\bv,\overline{\bw})$-mates.
}
\end{definition}
%%%%%
We clarify the notation $\lambda \not\equiv 0$. 
Throughout this paper, $\lambda \not\equiv 0$ means that $\{t \in I | \lambda(t) \not=0\}$ is a dense subset of $I$.
Then $\lambda$ is not identically zero for any non-trivial subintervals of $I$. 
It follows that $\gamma$ and $\overline{\gamma}$ are different space curves for any non-trivial subintervals of $I$. 
Note that if $\lambda$ is constant, then $\lambda \not\equiv 0$ means that $\lambda$ is a non-zero constant.
\par
We give all characterizations of Bertrand type curves of non-degenerate regular curves. 
Let $\gamma: I \to \R^3$ be a non-degenerate curve with curvature and torsion $(\kappa,\tau)$.
By a parameter change, we may assume that $s$ is the arc-length parameter of  $\gamma$. 
Note that $s$ is not the arc-length parameter of $\overline{\gamma}$.

%%%%%
\begin{proposition}\label{tt-Bertrand-type}
$\gamma$ is not a $(\bt,\overline{\bt})$-Bertrand type curve.
\end{proposition}
%%%%%
\demo Suppose that $\gamma$ is a $(\bt,\overline{\bt})$-Bertrand type curve. By differentiating $\overline{\gamma}(s)=\gamma(s)+\lambda(s)\bt(s)$, we have
$
|\dot{\overline{\gamma}}(s)|\overline{\bt}(s)=(1+\lambda'(s))\bt (s)+\lambda(s)\kappa(s)\bn(s).
$
Since $\bt(s)=\pm\overline{\bt}(s)$, we have $\lambda(s)\kappa(s)=0$ for all $s\in I$. Furthermore, since $\kappa(s)>0$, $\lambda(s)=0$ for all $s\in I$ and hence $\overline{\gamma}(t)=\gamma(t)$ for all $t\in I$. It follows that $\gamma$ is not a $(\bt,\overline{\bt})$-Bertrand type curve.
\enD
%%%%%

%%%%%
\begin{theorem}\label{tn-Bertrand-type}
$\gamma$ is a $(\bt,\overline{\bn})$-Bertrand type curve if and only if  $\tau(s)=0$ and there exists a constant $c \in \R$ such that $-s+c \not=0$ for all $s \in I$.
\end{theorem}
%%%%%
\demo
Suppose that $\gamma$ is a $(\bt,\overline{\bn})$-Bertrand type curve. By differentiating $\overline{\gamma}(s)=\gamma(s)+\lambda(s)\bt(s)$, we have
$
|\dot{\overline{\gamma}}(s)|\overline{\bt}(s)=(1+\lambda'(s))\bt (s)+\lambda(s)\kappa(s)\bn(s).
$
Since $\bt(s)=\pm\overline{\bn}(s)$, we have $1+\lambda'(s)=0$ for all $s\in I$. Therefore 
there exists a constant $c\in \R$ such that $\lambda(s)=-s+c$. 
If $\bt(s)=\overline{\bn}(s)$, there exists a smooth function $\theta:I \to \R$ such that
$$
\left(
\begin{array}{c}
\overline{\bb}(s)\\
\overline{\bt}(s)
\end{array}
\right)
=
\left(
\begin{array}{cc}
\cos\theta(s) & -\sin\theta(s)\\
\sin\theta(s) & \cos\theta(s)
\end{array}
\right)
\left(
\begin{array}{c}
\bn(s)\\
\bb(s)
\end{array}
\right).
$$
Then $|\dot{\overline{\gamma}}(s)|\cos\theta(s)=0$ and $|\dot{\overline{\gamma}}(s)|\sin\theta(s)=\lambda(s)\kappa(s)$. Since $|\dot{\overline{\gamma}}(s)|\neq0$, 
we have $\cos\theta(s)=0$, $\sin\theta(s)=\pm1$ and hence $\overline{\bb}(s)=\mp\bb(s)$ and $\overline{\bt}(s)=\pm\bn(s)$ for all $s\in I$. Furthermore, since $|\dot{\overline{\gamma}}(s)|>0$ and $\kappa(s)>0$, 
we have $\lambda(s)=-s+c\neq0$ for all $s\in I$. By differentiating $\overline{\bb}(s)=\mp\bb(s)$, we have 
$
-|\dot{\overline{\gamma}}(s)|\overline{\tau}(s)\overline{\bn}(s)=\pm\tau(s)\bn(s).
$
Since $\bt(s)=\overline{\bn}(s)$, we have $\tau(s)=0$ and $\overline{\tau}(s)=0$ for all $s\in I$. By differentiating $\overline{\bt}(s)=\mp\bn(s)$, we have
$
|\dot{\overline{\gamma}}(s)|\overline{\kappa}(s)\overline{\bn}(s)=\mp\kappa(s)\bt(s).
$
Since $\bt(s)=\overline{\bn}(s)$, we have $|\dot{\overline{\gamma}}(s)|\overline{\kappa}(s)=\mp\kappa(s)$ for all $s\in I$. Furthermore, since $|\dot{\overline{\gamma}}(s)|\overline{\kappa}(s)>0$, we have $\overline{\bt}(s)=-\bn(s)$, $\sin\theta(s)=-1$ and $\overline{\bb}(s)=\bb(s)$. Hence, we have
\begin{align*}
|\dot{\overline{\gamma}}(s)|=-\lambda(s)\kappa(s), \  
\kappa(s)=-\frac{|\dot{\overline{\gamma}}(s)|}{\lambda(s)}=\frac{|\dot{\overline{\gamma}}(s)|}{s-c}.
\end{align*}
Since $\kappa(s)>0$, if $\bt(s)=\overline{\bn}(s)$, we have $s-c>0$ for all $s\in I$.\par
On the other hand, if $\bt(s)=-\overline{\bn}(s)$, there exists a smooth function $\theta:I \to \R$ such that
$$
\left(
\begin{array}{c}
\overline{\bb}(s)\\
\overline{\bt}(s)
\end{array}
\right)
=
\left(
\begin{array}{cc}
\cos\theta(s) & -\sin\theta(s)\\
\sin\theta(s) & \cos\theta(s)
\end{array}
\right)
\left(
\begin{array}{c}
\bb(s)\\
\bn(s)
\end{array}
\right).
$$
Then $|\dot{\overline{\gamma}}(s)|\sin\theta(s)=0$ and $|\dot{\overline{\gamma}}(s)|\cos\theta(s)=\lambda(s)\kappa(s)$. Since $|\dot{\overline{\gamma}}(s)|\neq0$, 
we have $\sin\theta(s)=0$, $\cos\theta(s)=\pm1$ and hence $\overline{\bb}(s)=\pm\bb(s)$ and $\overline{\bt}(s)=\pm\bn(s)$ for all $s\in I$. Furthermore, since $|\dot{\overline{\gamma}}(s)|>0$ and $\kappa(s)>0$, we have $\lambda(s)=-s+c\neq0$ for all $s\in I$. By differentiating $\overline{\bb}(s)=\pm\bb(s)$, we have 
$
-|\dot{\overline{\gamma}}(s)|\overline{\tau}(s)\overline{\bn}(s)=\mp\tau(s)\bn(s).
$
Since $\bt(s)=-\overline{\bn}(s)$, we have $\tau(s)=0$ and $\overline{\tau}(s)=0$ for all $s\in I$. By differentiating $\overline{\bt}(s)=\pm\bn(s)$, we have
$
|\dot{\overline{\gamma}}(s)|\overline{\kappa}(s)\overline{\bn}(s)=\mp\kappa(s)\bt(s).
$
Since $\bt(s)=-\overline{\bn}(s)$, we have $|\dot{\overline{\gamma}}(s)|\overline{\kappa}(s)=\pm\kappa(s)$ for all $s\in I$. Furthermore, since $|\dot{\overline{\gamma}}(s)|\overline{\kappa}(s)>0$, we have $\overline{\bt}(s)=\bn(s)$, $\cos\theta(s)=1$ and $\overline{\bb}(s)=\bb(s)$. Hence, we have
\begin{align*}
|\dot{\overline{\gamma}}(s)|=\lambda(s)\kappa(s), \ 
\kappa(s)=\frac{|\dot{\overline{\gamma}}(s)|}{\lambda(s)}=\frac{|\dot{\overline{\gamma}}(s)|}{-s+c}.
\end{align*}
Since $\kappa(s)>0$, if $\bt(s)=-\overline{\bn}(s)$, we have $-s+c>0$ for all $s\in I$. It follows that we have $-s+c\neq0$ for all $s\in I$.\par
Conversely, suppose that $\tau(s)=0$ and there exists a constant $c\in\R$ such that $-s+c\neq0$ for all $s\in I$. Set $\overline{\gamma}(s)=\gamma(s)+(-s+c)\bt(s)$. By a direct calculation, we have
\begin{align*}
\dot{\overline{\gamma}}(s)&=|\dot{\overline{\gamma}}(s)|\overline{\bt}(s)=(-s+c)\kappa(s)\bn(s),\\
\ddot{\overline{\gamma}}(s)&=\left((-s+c)\kappa(s)\right)'\bn(s)+(-s+c)\kappa^2(s)\bt(s),\\
\dot{\overline{\gamma}}(s)&\times\ddot{\overline{\gamma}}(s)=(-s+c)^2\kappa^3(s)\bb(s)\neq0.
\end{align*}
It follows that $\overline{\gamma}$ is a non-degenerate curve. Since $|\dot{\overline{\gamma}}(s)|=\sqrt{(-s+c)^2\kappa^2(s)}$, we have
\begin{align*}
\overline{\bt}(s)&=\frac{(-s+c)\kappa(s)}{\sqrt{(-s+c)^2\kappa^2(s)}}\bn(s)=\frac{(-s+c)\kappa(s)}{|-s+c|\kappa(s)}\bn(s)=\pm\bn(s).
\end{align*}
By differentiating $\overline{\bt}(s)=\pm\bn(s)$, we have $|\dot{\overline{\gamma}}(s)|\overline{\kappa}(s)\overline{\bn}(s)=\mp\kappa(s)\bt(s)$. Hence, we have $\overline{\bn}(s)=\pm\bt(s)$. It follows that $\gamma$ is a $(\bt, \overline{\bn})$-Bertrand type curve.
\enD
%%%%%
\begin{remark} {\rm
If $\gamma$ is a $(\bt, \overline{\bn})$-Bertrand type curve, then $\overline{\gamma}$ is an involute of $\gamma$ (cf. \cite{Bruce-Giblin, Gibson,Gray}). Moreover, the curvature $\overline{\kappa}$ and the torsion $\overline{\tau}$ of $\overline{\gamma}$ are given by $\overline{\kappa}(s)=1/|-s+c|, \overline{\tau}(s)=\tau(s)=0$ for all $s\in I$.
}
\end{remark}

%%%%%
\begin{proposition}\label{tb-Bertrand-type}
$\gamma$ is not a $(\bt,\overline{\bb})$-Bertrand type curve.
\end{proposition}
%%%%%
\demo
Suppose that $\gamma$ is a $(\bt,\overline{\bb})$-Bertrand type curve. By differentiating $\overline{\gamma}(s)=\gamma(s)+\lambda(s)\bt(s)$, we have
$
|\dot{\overline{\gamma}}(s)|\overline{\bt}(s)=(1+\lambda'(s))\bt (s)+\lambda(s)\kappa(s)\bn(s).
$
If $\bt(s)=\overline{\bb}(s)$, there exists a smooth function $\theta:I \to \R$ such that
$$
\left(
\begin{array}{c}
\overline{\bt}(s)\\
\overline{\bn}(s)
\end{array}
\right)
=
\left(
\begin{array}{cc}
\cos\theta(s) & -\sin\theta(s)\\
\sin\theta(s) & \cos\theta(s)
\end{array}
\right)
\left(
\begin{array}{c}
\bn(s)\\
\bb(s)
\end{array}
\right).
$$
Then $|\dot{\overline{\gamma}}(s)|\sin\theta(s)=0$ and $|\dot{\overline{\gamma}}(s)|\cos\theta(s)=\lambda(s)\kappa(s)$. Since $|\dot{\overline{\gamma}}(s)|\neq0$, 
we have $\sin\theta(s)=0$, $\cos\theta(s)=\pm1$ and hence $\overline{\bt}(s)=\pm\bn(s)$ for all $s\in I$. By differentiating $\overline{\bt}(s)=\pm\bn(s)$, we have
$
 -|\dot{\overline{\gamma}}(s)|\overline{\kappa}(s)\overline{\bn}(s)=\mp\kappa(s)\bt(s)\pm\tau(s)\bb(s).
$
Since $\bt(s)=\overline{\bb}(s)$, we have $\kappa(s)=0$ for all $s\in I$. Therefore $\gamma$ is degenerate curve. It is a contradict of the definition.\par
On the other hand, if $\bt(s)=-\overline{\bb}(s)$, there exists a smooth function $\theta:I \to \R$ such that
$$
\left(
\begin{array}{c}
\overline{\bt}(s)\\
\overline{\bn}(s)
\end{array}
\right)
=
\left(
\begin{array}{cc}
\cos\theta(s) & -\sin\theta(s)\\
\sin\theta(s) & \cos\theta(s)
\end{array}
\right)
\left(
\begin{array}{c}
\bb(s)\\
\bn(s)
\end{array}
\right).
$$
Then $|\dot{\overline{\gamma}}(s)|\cos\theta(s)=0$ and $|\dot{\overline{\gamma}}(s)|\sin\theta(s)=-\lambda(s)\kappa(s)$. Since $|\dot{\overline{\gamma}}(s)|\neq0$, 
we have $\cos\theta(s)=0$, $\sin\theta(s)=\pm1$ and hence $\overline{\bt}(s)=\mp\bn(s)$ for all $s\in I$. By differentiating $\overline{\bt}(s)=\mp\bn(s)$, we have
$
-|\dot{\overline{\gamma}}(s)|\overline{\kappa}(s)\overline{\bn}(s)=\pm\kappa(s)\bt(s)\mp\tau(s)\bb(s).
$
Since $\bt(s)=-\overline{\bb}(s)$, we have $\kappa(s)=0$ for all $s\in I$. Therefore $\gamma$ is degenerate curve. It is a contradict of the definition. It follows that $\gamma$ is not a $(\bt,\overline{\bb})$-Bertrand type curve.
\enD
%%%%%

%%%%%
\begin{theorem}\label{nt-Bertrand-type}
$\gamma$ is a $(\bn,\overline{\bt})$-Bertrand type curve if and only if $\tau(s)=0$ and $\kappa'(s)\neq0$ for all $s \in I$.
\end{theorem}
%%%%%
\demo Suppose that $\gamma$ is a $(\bn,\overline{\bt})$-Bertrand type curve. By differentiating $\overline{\gamma}(s)=\gamma(s)+\lambda(s)\bn(s)$, we have
$
|\dot{\overline{\gamma}}(s)|\overline{\bt}(s)=(1-\lambda(s)\kappa(s))\bt(s)+\lambda'(s)\bn(s)+\lambda(s)\tau(s)\bb(s).
$
Since $\bn(s)=\pm\overline{\bt}(s)$, we have $|\dot{\overline{\gamma}}(s)|=\pm\lambda'(s)$, $1-\lambda(s)\kappa(s)=0$ and $\lambda(s)\tau(s)=0$ for all $s \in I$. Since $\kappa(s)>0$, we have $\lambda(s)=1/\kappa(s)>0$ and hence $\tau(s)=0$ for all $s\in I$. By differentiating $\lambda(s)=1/\kappa(s)$, we have $\lambda'(s)=-{\kappa'(s)}/{\kappa^2(s)}$. Since $|\dot{\overline{\gamma}}(s)|=\pm\lambda'(s)$ and $|\dot{\overline{\gamma}}(s)|>0$, we have 
$|\dot{\overline{\gamma}}(s)|={\kappa'(s)}/{\kappa^2(s)}\neq0$ for all $s\in I$. \par
Conversely, suppose that $\tau(s)=0$ and $\kappa'(s)\neq0$ for all $s\in I$. Set $\overline{\gamma}(s)=\gamma(s)+\bn(s)/\kappa(s)$. By a direct calculation, we have 
\begin{align*}
\dot{\overline{\gamma}}(s)&=|\dot{\overline{\gamma}}(s)|\overline{\bt}(s)=-\frac{\kappa'(s)}{\kappa^2(s)}\bn(s),\ 
\ddot{\overline{\gamma}}(s)=\left(-\frac{\kappa'(s)}{\kappa^2(s)}\right)'\bn(s)+\frac{\kappa'(s)}{\kappa(s)}\bt(s),\\
\dot{\overline{\gamma}}(s)&\times\ddot{\overline{\gamma}}(s)=\frac{\kappa'(s)^2}{\kappa^3(s)}\bb(s)\neq0.
\end{align*}
It follows that $\overline{\gamma}$ is a non-degenerate curve. Since $|\dot{\overline{\gamma}}(s)|=|\kappa'(s)|/\kappa^2(s)$, we have 
$$
\frac{|\kappa'(s)|}{\kappa^2(s)}\overline{\bt}(s)=-\frac{\kappa'(s)}{\kappa^2(s)}\bn(s).
$$
Hence, we have $\overline{\bt}(s)=\pm\bn(s)$. It follows that $\gamma$ is a $(\bn,\overline{\bt})$-Bertrand type curve.
\enD
%%%%%
\begin{remark} {\rm 
If $\gamma$ is a $(\bn, \overline{\bt})$-Bertrand type curve, then $\overline{\gamma}$ is an evolute of $\gamma$  (cf. \cite{Bruce-Giblin, Gibson,Gray}). 
Moreover, the curvature $\overline{\kappa}$ and the torsion $\overline{\tau}$ of $\overline{\gamma}$ are given by $\overline{\kappa}(s)=\kappa^3(s)/|\kappa'(s)|, \overline{\tau}(s)=\tau(s)=0$ for all $s \in I$.
}
\end{remark}

%%%%%
\begin{proposition}\label{bt-Bertrand-type}
$\gamma$ is not a $(\bb,\overline{\bt})$-Bertrand type curve.
\end{proposition}
%%%%%
\demo Suppose that $\gamma$ is a $(\bb,\overline{\bt})$-Bertrand type curve. By differentiating $\overline{\gamma}(s)=\gamma(s)+\lambda(s)\bb(s)$, we have
$
|\dot{\overline{\gamma}}(s)|\overline{\bt}(s)=\bt(s)+\lambda'(s)\bb(s)-\lambda(s)\tau(s)\bn(s).
$
Since $\bb(s)=\pm\overline{\bt}(s)$, we have $|\dot{\overline{\gamma}}(s)|=\pm\lambda'(s)$, $\bt(s)-\lambda(s)\tau(s)\bn(s)=0$ and hence $\bt(s)=\lambda(s)\tau(s) \bn(s)$ for all $s\in I$. It is a contradict of a moving frame of $\gamma(s)$.
 It follows that $\gamma$ is not a $(\bb,\overline{\bt})$-Bertrand type curve.
\enD
%%%%%

In \cite{Liu-Wang}, they investigate $(\bb,\overline{\bn})$-Bertrand type as Mannheim partner.
However, it seems to be missed the non-degenerate condition.
%%%%%
\begin{theorem}\label{bn-Bertrand-type}
$\gamma$ is a $(\bb,\overline{\bn})$-Bertrand type curve if and only if $\tau(s)\neq0$ and there exists a non-zero constant $A$ such that $A\tau'(s)=\kappa(s)(A^2\tau^2(s)+1)$ for all $s\in I$. 
\end{theorem}
%%%%%
\demo Suppose that $\gamma$ is a $(\bb,\overline{\bn})$-Bertrand type curve. By differentiating $\overline{\gamma}(s)=\gamma(s)+\lambda(s)\bb(s)$, we have 
$
|\dot{\overline{\gamma}}(s)|\overline{\bt}(s)=\bt(s)+\lambda'(s)\bb(s)-\lambda(s)\tau(s)\bn(s).
$
Since $\bb(s)=\pm\overline{\bn}(s)$, we have $\lambda'(s)=0$ for all $s \in I$. Therefore $\lambda$ is a constant. If $\lambda=0$, then $\overline{\gamma}(s)=\gamma(s)$ for all $s \in I$. 
Hence, $\lambda$ is a non-zero constant. We rewrite $\lambda$ as $A$. 
By differentiating $\overline{\gamma}(s)=\gamma(s)+A\bb(s)$, 
we have 
$
|\dot{\overline{\gamma}}(s)|\overline{\bt}(s)=\bt(s)-A\tau(s)\bn(s).
$
If $\bb(s)=\overline{\bn}(s)$, there exists a smooth function $\theta:I \to \R$ such that
$$
\left(
\begin{array}{c}
\overline{\bb}(s)\\
\overline{\bt}(s)
\end{array}
\right)
=
\left(
\begin{array}{cc}
\cos\theta(s) & -\sin\theta(s)\\
\sin\theta(s) & \cos\theta(s)
\end{array}
\right)
\left(
\begin{array}{c}
\bt(s)\\
\bn(s)
\end{array}
\right).
$$
Then $|\dot{\overline{\gamma}}(s)|\sin\theta(s)=1$ and $|\dot{\overline{\gamma}}(s)|\cos\theta(s)=-A\tau(s)$. 
It follows that $\sin\theta(s)>0$ and
\begin{eqnarray}\label{equation7}
\cos\theta(s)+A\tau(s)\sin\theta(s)=0  
\end{eqnarray}
for all $s\in I$. By differentiating $\overline{\bt}(s)=\sin\theta(s)\bt(s)+\cos\theta(s)\bn(s)$, we have
$
 |\dot{\overline{\gamma}}(s)|\overline{\kappa}(s)\overline{\bn}(s)=(\kappa(s)-\theta'(s))\overline{\bb}(s)+\tau(s)\cos\theta(s)\bb(s).
$
Since $\bb(s)=\overline{\bn}(s)$, we have $\kappa(s)-\theta'(s)=0$ and $|\dot{\overline{\gamma}}(s)|\overline{\kappa}(s)=\tau(s)\cos\theta(s)$ for all $s\in I$. Since $|\dot{\overline{\gamma}}(s)|\overline{\kappa}(s)>0$, we have $\tau(s)\cos\theta(s)>0$ for all $s\in I$. If $\tau(s)=0$ at a point $s\in I$, then $|\dot{\overline{\gamma}}(s)|\overline{\kappa}(s)=0$. It is a contradict of $|\dot{\overline{\gamma}}(s)|\overline{\kappa}(s)>0$. Hence, $\tau(s)\neq0$ for all $s\in I$. By differentiating \eqref{equation7}, we have
$$
\theta'(s)(A\tau(s)\cos\theta(s)-\sin\theta(s))+A\tau'(s)\sin\theta(s)=0.
$$
By \eqref{equation7}, we have $\kappa(s)(A^2\tau^2(s)+1)=A\tau'(s)$ for all $s \in I$. 
\par
On the other hand, if $\bb(s)=-\overline{\bn}(s)$, there exists a smooth function $\theta:I \to \R$ such that
$$
\left(
\begin{array}{c}
\overline{\bb}(s)\\
\overline{\bt}(s)
\end{array}
\right)
=
\left(
\begin{array}{cc}
\cos\theta(s) & -\sin\theta(s)\\
\sin\theta(s) & \cos\theta(s)
\end{array}
\right)
\left(
\begin{array}{c}
\bn(s)\\
\bt(s)
\end{array}
\right).
$$
Then $|\dot{\overline{\gamma}}(s)|\cos\theta(s)=1$ and $|\dot{\overline{\gamma}}(s)|\sin\theta(s)=-A\tau(s)$. It follows that $\cos\theta(s)>0$ and
\begin{eqnarray}\label{equation8}
\sin\theta(s)+A\tau(s)\cos\theta(s)=0  
\end{eqnarray}
for all $s\in I$. By differentiating $\overline{\bt}(s)=\sin\theta(s)\bn(s)+\cos\theta(s)\bt(s)$, we have
$
 |\dot{\overline{\gamma}}(s)|\overline{\kappa}(s)\overline{\bn}(s)=(\kappa(s)+\theta'(s))\overline{\bb}(s)+\tau(s)\sin\theta(s)\bb(s).
$
Since $\bb(s)=-\overline{\bn}(s)$, we have $\kappa(s)+\theta'(s)=0$ and $|\dot{\overline{\gamma}}(s)|\overline{\kappa}(s)=-\tau(s)\sin\theta(s)$ for all $s\in I$. Since $|\dot{\overline{\gamma}}(s)|\overline{\kappa}(s)>0$, we have $\tau(s)\sin\theta(s)<0$ for all $s\in I$. If $\tau(s)=0$ at a point $s\in I$, then $|\dot{\overline{\gamma}}(s)|\overline{\kappa}(s)=0$. It is a contradict of $|\dot{\overline{\gamma}}(s)|\overline{\kappa}(s)>0$. Hence, $\tau(s)\neq0$ for all $s\in I$. 
By differentiating \eqref{equation8}, we have
$$
\theta'(s)(\cos\theta(s)+A\tau(s)\sin\theta(s))+A\tau'(s)\cos\theta(s)=0.
$$
By \eqref{equation8}, we have $\kappa(s)(A^2\tau^2(s)+1)=A\tau'(s)$ for all $s \in I$. 
\par
Conversely, suppose that $\tau(s)\neq0$ and there exists a non-zero constant $A$ such that $A\tau'(s)=\kappa(s)(A^2\tau^2(s)+1)$ for all $s\in I$. Set $\overline{\gamma}(s)=\gamma(s)+A\bb(s)$. By a direct calculation, we have 
\begin{align*}
\dot{\overline{\gamma}}(s) &=\bt(s)-A\tau(s)\bn(s), \\
\ddot{\overline{\gamma}}(s) &=A\kappa(s)\tau(s)\bt(s)+(\kappa(s)-A\tau'(s))\bn(s)-A\tau^2(s)\bb(s),\\
\dot{\overline{\gamma}}(s)\times\ddot{\overline{\gamma}}(s)&=A\tau^3(s)\bt(s)+A\tau^2(s)\bn(s)+(\kappa(s)-A\tau'(s)+A^2\kappa(s)\tau^2(s))\bb(s)\\
  &=A\tau^2(s)(A\tau(s)\bt(s)+\bn(s))\neq0.
\end{align*}
It follows that $\overline{\gamma}$ is a non-degenerate curve. Furthermore, we have
\begin{align*}
&(\dot{\overline{\gamma}}(s)\times\ddot{\overline{\gamma}}(s))\times \dot{\overline{\gamma}}(s)=-A\tau^2(s)(A^2\tau^2(s)+1)\bb(s).
\end{align*}
Since $\overline{\bn}(s)=\overline{\bb}(s)\times\overline{\bt}(s)$, we have
$$
\overline{\bn}(s)=\frac{\dot{\overline{\gamma}}(s)\times\ddot{\overline{\gamma}}(s)}{|\dot{\overline{\gamma}}(s)\times\ddot{\overline{\gamma}}(s)|}\times\frac{\dot{\overline{\gamma}}(s)}{|\dot{\overline{\gamma}}(s)|}
=-\frac{A\tau^2(s)(A^2\tau^2(s)+1)}{|A|\tau^2(s)(A^2\tau^2(s)+1)}\bb(s)
=\pm\bb(s).
$$
It follows that $\gamma$ is a $(\bb,\overline{\bn})$-Bertrand type curve.
\enD
%%%%%
\begin{proposition}
Let $\gamma$ and $\overline{\gamma}:I\rightarrow\R^3$ be different non-degenerate curves. Under the same assumptions in Theorem \ref{bn-Bertrand-type}, suppose that $\gamma$ and $\overline{\gamma}$ are $(\bb,\overline{\bn})$-Bertrand mates with $\overline{\gamma}(s)=\gamma(s)+A\bb(s)$, $\tau(s)\neq0$ and $A\tau'(s)=\kappa(s)(A^2\tau^2(s)+1)$ for all $s\in I$, where $A$ is a non-zero constant. Then the curvature $\overline{\kappa}$ and the torsion $\overline{\tau}$ of $\overline{\gamma}$ are given by
$$
\overline{\kappa}(s)=\frac{|A|\tau^2(s)}{1+A^2\tau^2(s)}, \ \overline{\tau}(s)=\frac{\tau(s)}{1+A^2\tau^2(s)}.
$$
\end{proposition}
\demo
Since $\overline{\gamma}(s)=\gamma(s)+A\bb(s)$, we have 
$\dot{\overline{\gamma}}(s)=\bt(s)-A\tau(s)\bn(s). $
Therefore, we have
\begin{align*}
\ddot{\overline{\gamma}}(s)&=A\tau(s)(\kappa(s)\bt(s)-A\tau(s)\kappa(s)\bn(s)-\tau(s)\bb(s)),\\
\dddot{\overline{\gamma}}(s)&=A\left((\kappa(s)\tau(s))'+A\kappa^2(s)\tau^2(s)\right)\bt(s)\\
&\quad +A\tau(s)\left(-2\kappa(s)\tau'(s)+\kappa^2(s)-A\kappa'(s)\tau(s)+\tau^2(s)\right)\bn(s)\\
&\quad -A\tau(s)\left(2\tau'(s)+A\kappa(s)\tau^2(s)\right)\bb(s).
\end{align*}
Since
\begin{align*}
|\dot{\overline{\gamma}}(s)|&=\sqrt{1+A^2\tau^2(s)},\\
\dot{\overline{\gamma}}(s) \times \ddot{\overline{\gamma}}(s)&=A\tau^2(s)\left(\bn(s)+A\tau(s)\bt(s)\right),\\
|\dot{\overline{\gamma}}(s) \times \ddot{\overline{\gamma}}(s)|&=|A|\tau^2(s)\sqrt{1+A^2\tau^2(s)},\\
{\rm det}(\dot{\overline{\gamma}}(s), \ddot{\overline{\gamma}}(s), \dddot{\overline{\gamma}}(s))&=A^2\tau^5(s),
\end{align*}
we have the curvature and the torsion as
\begin{align*}
\overline{\kappa}(s)=\frac{|\dot{\overline{\gamma}}(s) \times \ddot{\overline{\gamma}}(s)|}{|\dot{\overline{\gamma}}(s)|^3}=\frac{|A|\tau^2(s)}{1+A^2\tau^2(s)}, \ 
\overline{\tau}(s)=\frac{{\rm det}(\dot{\overline{\gamma}}(s), \ddot{\overline{\gamma}}(s), \dddot{\overline{\gamma}}(s))}{|\dot{\overline{\gamma}}(s) \times \ddot{\overline{\gamma}}(s)|^2}=\frac{\tau(s)}{1+A^2\tau^2(s)}.
\end{align*}
\enD
%%%%%
\begin{theorem}\label{bb-Bertrand-type}
$\gamma$ is a $(\bb,\overline{\bb})$-Bertrand type curve if and only if $\tau(s)=0$ 
for all $s \in I$.
\end{theorem}
%%%%%
\demo Suppose that $\gamma$ is a $(\bb,\overline{\bb})$-Bertrand type curve. By differentiating $\overline{\gamma}(s)=\gamma(s)+\lambda(s)\bb(s)$, we have 
$
|\dot{\overline{\gamma}}(s)|\overline{\bt}(s)=\bt(s)+\lambda'(s)\bb(s)-\lambda(s)\tau(s)\bn(s).
$
Since $\bb(s)=\pm\overline{\bb}(s)$, we have $\lambda'(s)=0$ for all $s\in I$. Therefore $\lambda$ is a constant. If $\lambda(s)=0$, then $\overline{\gamma}(s)=\gamma(s)$ for all $s\in I$. 
Hence, $\lambda$ is a non-zero constant. 
We rewrite $\lambda$ as $A$. 
%Note that $s$ is not the arc-length parameter of $\overline{\gamma}$. 
By differentiating $\overline{\gamma}(s)=\gamma(s)+A\bb(s)$, we have
$
|\dot{\overline{\gamma}}(s)|\overline{\bt}(s)=\bt(s)-A\tau(s)\bn(s).
$
If $\bb(s)=\overline{\bb}(s)$, there exists a smooth function $\theta:I \to \R$ such that
$$
\left(
\begin{array}{c}
\overline{\bt}(s)\\
\overline{\bn}(s)
\end{array}
\right)
=
\left(
\begin{array}{cc}
\cos\theta(s) & -\sin\theta(s)\\
\sin\theta(s) & \cos\theta(s)
\end{array}
\right)
\left(
\begin{array}{c}
\bt(s)\\
\bn(s)
\end{array}
\right).
$$
Then $|\dot{\overline{\gamma}}(s)|\cos\theta(s)=1$ and $|\dot{\overline{\gamma}}(s)|\sin\theta(s)=A\tau(s)$. It follows that $\cos\theta(s)>0$ for all $s\in I$. By differentiating $\overline{\bt}(s)=\cos\theta(s)\bt(s)-\sin\theta(s)\bn(s)$, we have
$$
 |\dot{\overline{\gamma}}(s)|\overline{\kappa}(s)\overline{\bn}(s)=(\kappa(s)-\theta'(s))\overline{\bn}(s)-\tau(s)\sin\theta(s)\bb(s).
$$
Since $\bb(s)=\overline{\bb}(s)$, we have $\tau(s)\sin\theta(s)=0$. If $\tau(s)\neq0$ at a point $s\in I$, then $\sin\theta(s)=0$. It is contradict of $|\dot{\overline{\gamma}}(s)|\sin\theta(s)=A\tau(s)$. It follows that we have $\tau(s)=0$ for all $s\in I$. \par
On the other hand, if $\bb(s)=-\overline{\bb}(s)$, there exists a smooth function $\theta:I \to \R$ such that
$$
\left(
\begin{array}{c}
\overline{\bt}(s)\\
\overline{\bn}(s)
\end{array}
\right)
=
\left(
\begin{array}{cc}
\cos\theta(s) & -\sin\theta(s)\\
\sin\theta(s) & \cos\theta(s)
\end{array}
\right)
\left(
\begin{array}{c}
\bn(s)\\
\bt(s)
\end{array}
\right).
$$
Then $|\dot{\overline{\gamma}}(s)|\sin\theta(s)=-1$ and $|\dot{\overline{\gamma}}(s)|\cos\theta(s)=-A\tau(s)$. It follows that $\sin\theta(s)<0$ for all $s\in I$. By differentiating $\overline{\bt}(s)=\cos\theta(s)\bn(s)-\sin\theta(s)\bt(s)$, we have
$$
 |\dot{\overline{\gamma}}(s)|\overline{\kappa}(s)\overline{\bn}(s)=-(\kappa(s)+\theta'(s))\overline{\bn}(s)+\tau(s)\cos\theta(s)\bb(s).
$$
Since $\bb(s)=-\overline{\bb}(s)$, we have $\tau(s)\cos\theta(s)=0$. If $\tau(s)\neq0$ at a point $s\in I$, then $\cos\theta(s)=0$. It is contradict of $|\dot{\overline{\gamma}}(s)|\cos\theta(s)=-A\tau(s)$. It follows that $\tau(s)=0$ for all $s\in I$ and $\cos\theta(s)=0$ and hence $\sin\theta(s)=-1$. Hence, we have $\overline{\bt}(s)=\bt(s)$ and $\overline{\bn}(s)=-\bn(s)$. By differentiating $\overline{\bt}(s)=\bt(s)$, we have
$
|\dot{\overline{\gamma}}(s)|\overline{\kappa}(s)\overline{\bn}(s)=\kappa(s)\bn(s).
$
Since $|\dot{\overline{\gamma}}(s)|>0$, $\overline{\kappa}(s)>0$ and $\kappa(s)>0$, we have $\overline{\bn}(s)=\bn(s)$. It is contradict of $\overline{\bn}(s)=-\bn(s)$. 
Hence, this case does not occur.
%It follows that if $\overline{\bb}(s)=\bb(s)$, $\gamma$ is a $(\bb,\overline{\bb})$-Bertrand type curve. 
\par
Conversely, suppose that $\tau(s)=0$ for all $s\in I$. Set $\overline{\gamma}(s)=\gamma(s)+A\bb(s)$, where $A$ is a non-zero constant. By a direct calculation, we have
\begin{align*}
\dot{\overline{\gamma}}(s)=\bt(s), \ \ddot{\overline{\gamma}}(s)=\kappa(s)\bn(s), \ 
\dot{\overline{\gamma}}(s)\times\ddot{\overline{\gamma}}(s)=\kappa(s)\bb(s)\neq0.
\end{align*}
It follows that $\overline{\gamma}$ is a non-degenerate curve. Since $\overline{\bb}(s)=\dot{\overline{\gamma}}(s)\times\ddot{\overline{\gamma}}(s)/|\dot{\overline{\gamma}}(s)\times\ddot{\overline{\gamma}}(s)|$, we have
$$
\overline{\bb}(s)=\frac{\dot{\overline{\gamma}}(s)\times\ddot{\overline{\gamma}}(s)}{|\dot{\overline{\gamma}}(s)\times\ddot{\overline{\gamma}}(s)|}=\frac{\kappa(s)}{\kappa(s)}\bb(s)=\bb(s)
$$
It follows that $\gamma$ is a $(\bb,\overline{\bb})$-Bertrand type curve.
\enD
%%%%%
\begin{remark} {\rm 
If $\gamma$ is a $(\bb, \overline{\bb})$-Bertrand type curve, $\overline{\gamma}$ is a parallel translation of $\gamma$. Moreover, the curvature $\overline{\kappa}$ and the torsion $\overline{\tau}$ of $\overline{\gamma}$ are given by $\overline{\kappa}(s)=\kappa(s), \overline{\tau}(s)=\tau(s)=0$ for all $s \in I$.
}
\end{remark}

%%%%%%%%%%%%%%%%% Section 4 %%%%%%%%%%%%%%%%%
\section{Bertrand framed curves}

Let $(\gamma,\nu_1,\nu_2)$ and $(\overline{\gamma},\overline{\nu}_1,\overline{\nu}_2):I \to \R^3 \times \Delta$ be framed curves. 
%%%%%
\begin{definition}\label{Bertrand-type-framed}{\rm
We say that $(\gamma,\nu_1,\nu_2)$ and $(\overline{\gamma},\overline{\nu}_1,\overline{\nu}_2)$ are {\it $(\bv,\overline{\bw})$-mates} if there exists a smooth function $\lambda:I \to \R$ with $\lambda \not\equiv 0$ such that $\overline{\gamma}(t)=\gamma(t)+\lambda(t)\bv(t)$ and $\bv(t)= \overline{\bw}(t)$ for all $t \in I$, 
where $\bv$ and $\bw$ are $\nu_1, \nu_2$ or $\mu$. 
We also say that  $(\gamma,\nu_1,\nu_2)$ is a {\it $(\bv,\overline{\bw})$-Bertrand framed curve} (or, {\it $(\bv,\overline{\bw})$-Bertrand-Mannheim framed curve}) if there exists another framed curve $(\overline{\gamma},\overline{\nu}_1,\overline{\nu}_2)$ such that $(\gamma,\nu_1,\nu_2)$ and $(\overline{\gamma},\overline{\nu}_1,\overline{\nu}_2)$ are $(\bv,\overline{\bw})$-mates.
}
\end{definition}
%%%%%
\par
A little bit the difference of the sign of $\overline{\bw}$ between the Definitions \ref{Bertrand-type-regular} and \ref{Bertrand-type-framed} comes from a flexibility of framed curves.
That is, if $(\gamma,\nu_1,\nu_2)$ is a framed curve, then $(\gamma,-\nu_1,\nu_2)$ and $(\gamma,\nu_1,-\nu_2)$ are also framed curves. 
Therefore, we may consider $\overline{\bw}$ up to sign.
\par
Let $(\gamma,\nu_1,\nu_2):I \to \R^3 \times \Delta$ be a framed curve with curvature $(\ell,m,n,\alpha)$. 
We give all characterizations of Bertrand framed curves.

%%%%%
\begin{theorem}\label{mu-mu}
$(\gamma,\nu_1,\nu_2):I \to \R^3 \times \Delta$ is a $(\mu,\overline{\mu})$-Bertrand framed curve if and only if $m(t)=n(t)=0$ for all $t \in I$.
\end{theorem}
%%%%%
\demo
Suppose that $(\gamma,\nu_1,\nu_2):I \to \R^3 \times \Delta$ is a $(\mu,\overline{\mu})$-Bertrand framed curve. 
Then there exists a framed curve $(\overline{\gamma},\overline{\nu}_1,\overline{\nu}_2)$ and a smooth function $\lambda:I \to \R$ with $\lambda \not\equiv 0$ such that 
$\overline{\gamma}(t)=\gamma(t)+\lambda(t)\mu(t)$ and $\mu(t)=\overline{\mu}(t)$ for all $t \in I$. 
By differentiating, we have $\overline{\alpha}(t) \overline{\mu}(t)=(\alpha(t)+\dot{\lambda}(t))\mu(t)-\lambda(t)m(t)\nu_1(t)-\lambda(t)n(t)\nu_2(t)=0$ for all $t \in I$. 
Since $\mu(t)=\overline{\mu}(t)$, we have $\overline{\alpha}(t)=\alpha(t)+\dot{\lambda}(t), \lambda(t)m(t)=0$ and $\lambda(t)n(t)=0$ for all $t \in I$. 
Since $\lambda \not\equiv 0$ and the continuous condition, we have $m(t)=n(t)=0$ for all $t \in I$. 
\par
Conversely, suppose that $(\gamma,\nu_1,\nu_2)$ is a framed curve with $m(t)=n(t)=0$ for all $t \in I$. 
Then $\mu$ is a constant vector $\bv$.
Consider a smooth map $(\overline{\gamma},\nu_1,\nu_2):I \to \R^3 \times \Delta$ with $\overline{\gamma}(t)=\gamma(t)+\lambda(t)\bv$. 
Since $\dot{\overline{\gamma}}(t)=(\alpha(t)+\dot{\lambda}(t))\bv$, we have $\dot{\overline{\gamma}}(t) \cdot \nu_1(t)=\dot{\overline{\gamma}}(t) \cdot \nu_2(t)=0$ for all $t \in I$. 
It follows that $(\overline{\gamma},\nu_1,\nu_2)$ is a framed curve and $\overline{\mu}(t)=\mu(t)=\bv$. 
Therefore, $(\gamma,\nu_1,\nu_2):I \to \R^3 \times \Delta$ is a $(\mu,\overline{\mu})$-Bertrand framed curve.
\enD
%%%%%
\begin{remark}{\rm 
If $(\gamma,\nu_1,\nu_2):I \to \R^3 \times \Delta$ is a $(\mu,\overline{\mu})$-Bertrand framed curve, then $\mu$ is a constant vector $\bv$ and $\gamma(t)=(\int \alpha(t) dt) \bv+\bc$, where $\bc \in \R^3$ is a constant vector. 
Therefore, $\gamma$ is a part of a line. 
}
\end{remark}
%%%%%

%%%%%
\begin{proposition}\label{mu-mu_curvature}
Suppose that $(\gamma,\nu_1,\nu_2)$ and $(\overline{\gamma},\overline{\nu}_1,\overline{\nu}_2): I \to \R^3 \times \Delta$ are {$(\mu,\overline{\mu})$-mates}, where $\overline{\gamma}(t)=\gamma(t)+\lambda(t)\mu(t),  \overline{\nu}_1(t)=\cos \theta(t){\nu_1}(t)-\sin \theta(t){\nu_2}(t), \overline{\nu}_2(t)=\sin \theta(t){\nu_1}(t)+\cos \theta(t){\nu_2}(t)$ and $\theta:I\rightarrow\R$ is a smooth function. 
Then the curvature $(\overline{\ell}, \overline{m}, \overline{n}, \overline{\alpha})$ of $(\overline{\gamma},\overline{\nu}_1,\overline{\nu}_2)$ is given by
$
\overline{\ell}(t)=\ell(t)-\dot{\theta}(t), 
\overline{m}(t)=0, 
\overline{n}(t)=0, 
\overline{\alpha}(t)=\alpha(t)+\dot{\lambda}(t).
$
\end{proposition}
%%%%%
\demo
%Since $\mu(t)=\overline{\mu}(t)$, there exists a smooth function $\theta:I \to \R$ such that $$\begin{pmatrix}\overline{\nu}_1(t) \\\overline{\nu}_2(t)\end{pmatrix}=\begin{pmatrix}\cos \theta(t) & -\sin \theta(t) \\\sin \theta(t) & \cos \theta(t)\end{pmatrix}\begin{pmatrix}{\nu}_1(t) \\{\nu}_2(t)\end{pmatrix}.$$Then we have $\overline{\nu}_1(t)=\cos \theta(t){\nu_1}(t)-\sin \theta(t){\nu_2}(t)$ and $\overline{\nu}_2(t)=\sin \theta(t){\nu_1}(t)+\cos \theta(t){\nu_2}(t)$. 
By differentiating $\overline{\nu}_1(t)=\cos \theta(t){\nu_1}(t)-\sin \theta(t){\nu_2}(t)$, we have 
$
\overline{\ell}(t)\overline{\nu}_2(t)+\overline{m}(t)\overline{\mu}(t)=(\ell(t)-\dot{\theta}(t))\overline{\nu}_2(t).
$
Then we have $\overline{\ell}(t)=\ell(t)-\dot{\theta}(t)$. Since $\nu_1(t)=\overline{\mu}(t)$, we have $\overline{m}(t)=0$. Moreover, by differentiating $\overline{\nu}_2(t)=\sin \theta(t){\nu_1}(t)+\cos \theta(t){\nu_2}(t)$, we have 
$
-\overline{\ell}(t)\overline{\nu}_1(t)+\overline{n}(t)\overline{\mu}(t)=(\ell(t)-\dot{\theta}(t))\overline{\nu}_1(t).
$
Since $\nu_1(t)=\overline{\mu}(t)$, we have $\overline{n}(t)=0$. By $\overline{\alpha}(t) \overline{\mu}(t)=(\alpha(t)+\dot{\lambda}(t))\mu(t)$ and $\mu(t)=\overline{\mu}(t)$, we have $\overline{\alpha}(t)=\alpha(t)+\dot{\lambda}(t)$.
\enD
%%%%%

%%%%%
\begin{theorem}\label{mu-nu1}
$(\gamma,\nu_1,\nu_2):I \to \R^3 \times \Delta$ is a $(\mu,\overline{\nu}_1)$-Bertrand framed curve if and only if there exists a smooth function $\theta:I \to \R$ such that 
$
m(t) \cos \theta(t)-n(t)\sin \theta(t)=0
$ 
for all $t \in I$ and $\int\alpha(t)dt\not\equiv 0$.
\end{theorem}
%%%%%
\demo
Suppose that $(\gamma,\nu_1,\nu_2):I \to \R^3 \times \Delta$ is a $(\mu,\overline{\nu}_1)$-Bertrand framed curve. 
Then there exists a framed curve $(\overline{\gamma},\overline{\nu}_1,\overline{\nu}_2)$ and a smooth function $\lambda:I \to \R$ with $\lambda \not\equiv 0$ such that 
$\overline{\gamma}(t)=\gamma(t)+\lambda(t)\mu(t)$ and $\mu(t)=\overline{\nu}_1(t)$ for all $t \in I$.
By differentiating, we have $\overline{\alpha}(t) \overline{\mu}(t)=(\alpha(t)+\dot{\lambda}(t))\mu(t)-\lambda(t)m(t)\nu_1(t)-\lambda(t)n(t)\nu_2(t)=0$ for all $t \in I$. 
Since $\mu(t)=\overline{\nu}_1(t)$, we have $\alpha(t)+\dot{\lambda}(t)=0$ for all $t \in I$. 
Moreover, there exists a smooth function $\theta:I \to \R$ such that 
$$
\begin{pmatrix}
\overline{\nu}_2(t) \\
\overline{\mu}(t)
\end{pmatrix}
=
\begin{pmatrix}
\cos \theta(t) & -\sin \theta(t) \\
\sin \theta(t) & \cos \theta(t)
\end{pmatrix}
\begin{pmatrix}
{\nu}_1(t) \\
{\nu}_2(t)
\end{pmatrix}.
$$
Then we have $\overline{\alpha}(t)\sin \theta(t)=-\lambda(t)m(t)$ and $\overline{\alpha}(t)\cos(t)=-\lambda(t)n(t)$.
It follows that $$\lambda(t)(m(t)\cos \theta(t)-n(t)\sin \theta(t))=0$$ for all $t \in I$. 
Since $\lambda \not\equiv 0$ and the continuous condition, we have $m(t)\cos \theta(t)-n(t)\sin \theta(t)=0$ for all $t \in I$.
\par
Conversely, suppose that there exists a smooth function $\theta:I \to \R$ such that $m(t)\cos \theta(t)-n(t)\sin \theta(t)=0$ for all $t \in I$ and $\int\alpha(t)dt\not\equiv 0$. 
Let $(\overline{\gamma},\overline{\nu}_1,\overline{\nu}_2):I \to \R^3 \times \Delta$ be 
$$
\overline{\gamma}(t)=\gamma(t)-\left(\int \alpha(t) dt\right)\mu(t),\ \overline{\nu}_1(t)=\mu(t),\ \overline{\nu}_2(t)=\cos \theta(t)\nu_1(t)-\sin \theta(t)\nu_2(t).
$$
Since $\dot{\overline{\gamma}}(t)=(\int \alpha(t) dt)(m(t)\nu_1(t)+n(t)\nu_2(t))$, we have  $\dot{\overline{\gamma}}(t) \cdot \overline{\nu}_1(t)=\dot{\overline{\gamma}}(t) \cdot \overline{\nu}_2(t)=0$ for all $t \in I$. 
It follows that $(\overline{\gamma},\overline{\nu}_1,\overline{\nu}_2)$ is a framed curve.
Therefore, $(\gamma,\nu_1,\nu_2)$ is a $(\mu,\overline{\nu}_1)$-Bertrand framed curve.
\enD
%%%%%

%%%%%
\begin{remark} {\rm
If $(\gamma,\nu_1,\nu_2)$ is a $(\mu, \overline{\nu}_1)$-Bertrand framed curve, then $\overline{\gamma}$ is an involute of the framed curve $(\gamma,\nu_1,\nu_2)$ (cf. \cite{Honda-Takahashi-Preprint}). 
}
\end{remark}
%%%%%

%%%%%
\begin{proposition}\label{mu-nu1_curvature}
Suppose that $(\gamma,\nu_1,\nu_2)$ and $(\overline{\gamma},\overline{\nu}_1,\overline{\nu}_2): I \to \R^3 \times \Delta$ are {$(\mu,\overline{\nu}_1)$-mates}, where $\overline{\gamma}(t)=\gamma(t)-(\int \alpha(t) dt)\mu(t),  \overline{\nu}_2(t)=\cos \theta(t){\nu_1}(t)-\sin \theta(t){\nu_2}(t)$ and $\overline{\mu}(t)=\sin \theta(t){\nu_1}(t)+\cos \theta(t){\nu_2}(t)$ and $\theta:I \to \R$ is a smooth function. 
Then the curvature $(\overline{\ell}, \overline{m}, \overline{n}, \overline{\alpha})$ of $(\overline{\gamma},\overline{\nu}_1,\overline{\nu}_2)$ is given by
\begin{align*}
  \overline{\ell}(t)&=-m(t)\cos\theta(t)+n(t)\sin\theta(t)=0, \\
  \overline{m}(t)&=-m(t)\sin\theta(t)-n(t)\cos\theta(t), \\
  \overline{n}(t)&=\ell(t)-\dot{\theta}(t), \\
  \overline{\alpha}(t)&=-\left(\int \alpha(t) dt\right)(m(t)\sin\theta(t)+n(t)\cos\theta(t)).
\end{align*}
\end{proposition}
%%%%%
\demo
By differentiating $\overline{\nu}_2(t)=\cos \theta(t){\nu_1}(t)-\sin \theta(t){\nu_2}(t)$, we have 
\begin{align*}
-\overline{\ell}(t)\overline{\nu}_1(t)+\overline{n}(t)\overline{\mu}(t)&=(\ell(t)-\dot{\theta}(t))\overline{\mu}(t)+(m(t)\cos\theta(t)-n(t)\sin\theta(t))\mu(t).
\end{align*}
Then we have $\overline{n}(t)=\ell(t)-\dot{\theta}(t)$. Since $\mu(t)=\overline{\nu}_1(t)$, we have $\overline{\ell}(t)=-m(t)\cos\theta(t)+n(t)\sin\theta(t)=0$. Moreover, by differentiating $\overline{\mu}(t)=\sin \theta(t){\nu_1}(t)+\cos \theta(t){\nu_2}(t)$, we have 
\begin{align*}
-\overline{m}(t)\overline{\nu}_1(t)-\overline{n}(t)\overline{\nu}_2(t)&=(\dot{\theta}(t)-\ell(t))\overline{\nu}_2(t)+(m(t)\sin\theta(t)+n(t)\cos\theta(t))\mu(t).
\end{align*}
Since $\nu_1(t)=\overline{\mu}(t)$, we have $\overline{m}(t)=-m(t)\sin\theta(t)-n(t)\cos\theta(t)$. By $\overline{\alpha}(t)\sin \theta(t)=-\lambda(t)m(t)$ and $\overline{\alpha}(t)\cos(t)=-\lambda(t)n(t)$, we have $\overline{\alpha}(t)=-\lambda(t)(m(t)\sin\theta(t)+n(t)\cos\theta(t))$.
\enD
%%%%%

%%%%%
\begin{theorem}\label{mu-nu2}
$(\gamma,\nu_1,\nu_2):I \to \R^3 \times \Delta$ is a $(\mu,\overline{\nu}_2)$-Bertrand framed curve if and only if $(\gamma,\nu_1,\nu_2):I \to \R^3 \times \Delta$ is a $(\mu,\overline{\nu}_1)$-Bertrand framed curve. 
\end{theorem}
%%%%%
\demo
If $(\gamma,\nu_1,\nu_2):I \to \R^3 \times \Delta$ is a $(\mu,\overline{\nu}_2)$-Bertrand framed curve, then there exists a framed curve $(\overline{\gamma},\overline{\nu}_1,\overline{\nu}_2)$ and a smooth function $\lambda:I \to \R$ with $\lambda \not\equiv 0$ such that 
$\overline{\gamma}(t)=\gamma(t)+\lambda(t)\mu(t)$ and $\mu(t)=\overline{\nu}_2(t)$ for all $t \in I$
Since $(\overline{\gamma},\overline{\nu}_2,\overline{\nu}_1)$ is also a framed curve, $(\gamma,\nu_1,\nu_2)$ is a $(\mu,\overline{\nu}_1)$-Bertrand framed curve and vice versa (see, Theorem \ref{relation} and Remark \ref{relation-remark}).
\enD
%%%%%
%%%%%
\begin{remark}\label{mu-nu1_mu-nu2}{\rm 
Let $(\gamma,\nu_1,\nu_2)$ and $(\overline{\gamma},\overline{\nu}_1,\overline{\nu}_2)$ be $(\mu,\overline{\nu}_1)$-mates. 
The curvature of $(\overline{\gamma},\overline{\nu}_1,\overline{\nu}_2)$ is given by $(\overline{\ell}, \overline{m}, \overline{n}, \overline{\alpha})$. 
Then $(\gamma,\nu_1,\nu_2)$ and $(\overline{\gamma},\overline{\nu}_2,\overline{\nu}_1)$ are $(\mu,\overline{\nu}_2)$-mates and the curvature of $(\overline{\gamma},\overline{\nu}_2,\overline{\nu}_1)$ is given by $(0, -\overline{n}, -\overline{m}, -\overline{\alpha})$ by Propositions \ref{change-frame} and \ref{mu-nu1_curvature}.
}
\end{remark}
%%%%%

%%%%%
\begin{theorem}\label{nu1-mu}
$(\gamma,\nu_1,\nu_2):I \to \R^3 \times \Delta$ is a $(\nu_1,\overline{\mu})$-Bertrand framed curve if and only if $\ell(t)=0$ and there exists a smooth function $\lambda : I\rightarrow\R$ with $\lambda \not\equiv 0$ such that $\alpha(t)+\lambda(t)m(t)=0$ for all $t\in I$. 
\end{theorem}
%%%%%
\demo
Suppose that $(\gamma,\nu_1,\nu_2):I \to \R^3 \times \Delta$ is a $(\nu_1,\overline{\mu})$-Bertrand framed curve. 
Then there exists a framed curve $(\overline{\gamma},\overline{\nu}_1,\overline{\nu}_2)$ and a smooth function $\lambda:I \to \R$ with $\lambda \not\equiv 0$ such that 
$\overline{\gamma}(t)=\gamma(t)+\lambda(t)\nu_1(t)$ and $\nu_1(t)=\overline{\mu}(t)$ for all $t \in I$.
By differentiating, we have $\overline{\alpha}(t) \overline{\mu}(t)=(\alpha(t)+\lambda(t)m(t))\mu(t)+\dot{\lambda}(t)\nu_1(t)+\lambda(t)\ell(t)\nu_2(t)=0$ for all $t \in I$. 
Since $\nu_1(t)=\overline{\mu}(t)$, we have $\alpha(t)+\lambda(t)m(t)=0$ and $\lambda(t)\ell(t)=0$ for all $t \in I$. 
Since $\lambda \not\equiv 0$, we have $\ell(t)=0$ for all $t \in I$. 
\par
Conversely, suppose that $\ell(t)=0$ and there exists a smooth function $\lambda : I\rightarrow\R$ such that $\alpha(t)+\lambda(t)m(t)=0$ for all $t\in I$. 
Let $(\overline{\gamma},\overline{\nu}_1,\overline{\nu}_2):I \to \R^3 \times \Delta$ be 
$
\overline{\gamma}(t)=\gamma(t)+\lambda(t)\nu_1(t),\ \overline{\nu}_1(t)=\nu_2(t),\ \overline{\nu}_2(t)=\mu(t).
$
Since $\dot{\overline{\gamma}}(t)=\dot{\lambda}(t)\nu_1(t)$, we have  $\dot{\overline{\gamma}}(t) \cdot \overline{\nu}_1(t)=\dot{\overline{\gamma}}(t) \cdot \overline{\nu}_2(t)=0$ for all $t \in I$. 
It follows that $(\overline{\gamma},\overline{\nu}_1,\overline{\nu}_2)$ is a framed curve. 
Moreover, since $\overline{\mu}(t)=\overline{\nu}_1(t)\times\overline{\nu}_2(t)$, we have $\overline{\mu}(t)={\nu}_1(t)$ for all $t\in I$. 
Therefore, $(\gamma,\nu_1,\nu_2)$ is a $(\nu_1,\overline{\mu})$-Bertrand framed curve.
\enD
%%%%%
%%%%%
\begin{remark} {\rm
If $(\gamma,\nu_1,\nu_2)$ is a $(\nu_1, \overline{\mu})$-Bertrand framed curve and $m(t) \not=0$ for all $t \in I$, then $\overline{\gamma}$ is a circular evolute with respect to $\nu_1$ of the framed curve $(\gamma,\nu_1,\nu_2)$ (cf. \cite{Honda-Takahashi-Preprint}). 
}
\end{remark}
%%%%%
%%%%%
\begin{proposition}\label{nu1-mu_curvature}
Suppose that $(\gamma,\nu_1,\nu_2)$ and $(\overline{\gamma},\overline{\nu}_1,\overline{\nu}_2): I \to \R^3 \times \Delta$ are {$(\nu_1,\overline{\mu})$-mates}, where $\overline{\gamma}(t)=\gamma(t)+\lambda(t)\nu_1(t), \overline{\nu}_1(t)=\cos \theta(t){\nu_2}(t)-\sin \theta(t){\mu}(t), \overline{\nu}_2(t)=\sin \theta(t){\nu_2}(t)+\cos \theta(t){\mu}(t)$ and $\theta:I \to \R$ is a smooth function. 
Then the curvature $(\overline{\ell}, \overline{m}, \overline{n}, \overline{\alpha})$ of $(\overline{\gamma},\overline{\nu}_1,\overline{\nu}_2)$ is given by
\begin{align*}
  \overline{\ell}(t)=n(t)-\dot{\theta}(t), \
  \overline{m}(t)=m(t)\sin\theta(t), \
  \overline{n}(t)=-m(t)\cos\theta(t), \
  \overline{\alpha}(t)=\dot{\lambda}(t).
\end{align*}
\end{proposition}
%%%%%
\demo
By differentiating $\overline{\nu}_1(t)=\cos \theta(t){\nu_2}(t)-\sin \theta(t){\mu}(t)$, we have 
\begin{align*}
\overline{\ell}(t)\overline{\nu}_2(t)+\overline{m}(t)\overline{\mu}(t)&=(n(t)-\dot{\theta}(t))\overline{\nu}_2(t)+m(t)\sin\theta(t)\nu_1(t).
\end{align*}
Then we have $\overline{\ell}(t)=n(t)-\dot{\theta}(t)$. Since $\nu_1(t)=\overline{\mu}(t)$, we have $\overline{m}(t)=m(t)\sin\theta(t)$. Moreover, by differentiating $\overline{\nu}_2(t)=\sin \theta(t){\nu_2}(t)+\cos \theta(t){\mu}(t)$, we have 
\begin{align*}
-\overline{\ell}(t)\overline{\nu}_1(t)+\overline{n}(t)\overline{\mu}(t)&=(\dot{\theta}(t)-n(t))\overline{\nu}_1(t)-m(t)\cos\theta(t)\nu_1(t).
\end{align*}
Since $\nu_1(t)=\overline{\mu}(t)$, we have $\overline{n}(t)=-m(t)\cos\theta(t)$. By $\overline{\alpha}(t) \overline{\mu}(t)=(\alpha(t)+\lambda(t)m(t))\mu(t)+\dot{\lambda}(t)\nu_1(t)+\lambda(t)\ell(t)\nu_2(t)=0$ and $\nu_1(t)=\overline{\mu}(t)$, we have $\overline{\alpha}(t)=\dot{\lambda}(t)$.
\enD
%%%%%

%%%%%
\begin{theorem}\label{nu2-nu1}
$(\gamma,\nu_1,\nu_2):I \to \R^3 \times \Delta$ is a $(\nu_2,\overline{\nu}_1)$-Bertrand framed curve if and only if there exists a non-zero constant $\lambda$ and a smooth function $\theta:I \to \R$ such that 
\begin{eqnarray*}
\lambda \ell(t) \sin \theta(t)+(\alpha(t)+\lambda n(t))\cos \theta(t)=0
\end{eqnarray*} 
for all $t \in I$.
\end{theorem}
%%%%%
\demo
Suppose that $(\gamma,\nu_1,\nu_2):I \to \R^3 \times \Delta$ is a $(\nu_2,\overline{\nu}_1)$-Bertrand framed curve. 
Then there exists a framed curve $(\overline{\gamma},\overline{\nu}_1,\overline{\nu}_2)$ and a smooth function $\lambda:I \to \R$ with $\lambda \not\equiv 0$ such that 
$\overline{\gamma}(t)=\gamma(t)+\lambda(t)\nu_2(t)$ and $\nu_2(t)=\overline{\nu}_1(t)$ for all $t \in I$.
By differentiating, we have $\overline{\alpha}(t) \overline{\mu}(t)=(\alpha(t)+\lambda(t)n(t))\mu(t)-\lambda(t)\ell(t)\nu_1(t)+\dot{\lambda}(t)\nu_2(t)=0$ for all $t \in I$. 
Since $\nu_2(t)=\overline{\nu}_1(t)$, 
we have $\dot{\lambda}(t)=0$ for all $t\in I$. Therefore $\lambda$ is a constant. If $\lambda=0$, then $\overline{\gamma}(t)=\gamma(t)$ for all $t\in I$. Hence, $\lambda$ is a non-zero constant.  
Moreover, there exists a smooth function $\theta:I \to \R$ such that 
$$
\begin{pmatrix}
\overline{\nu}_2(t) \\
\overline{\mu}(t)
\end{pmatrix}
=
\begin{pmatrix}
\cos \theta(t) & -\sin \theta(t) \\
\sin \theta(t) & \cos \theta(t)
\end{pmatrix}
\begin{pmatrix}
{\nu}_1(t) \\
{\nu}_2(t)
\end{pmatrix}.
$$
Then we have $\overline{\alpha}(t)\sin \theta(t)=\alpha(t)+\lambda n(t)$ and $\overline{\alpha}(t)\cos(t)=-\lambda\ell(t)$.
It follows that $$\lambda\ell(t)\sin \theta(t)+(\alpha(t)+\lambda n(t))\cos \theta(t)=0$$ for all $t \in I$. 
\par
Conversely, suppose that there exists a non-zero constant $\lambda$ and a smooth function $\theta:I \to \R$ such that 
$\lambda \ell(t) \sin \theta(t)+(\alpha(t)+\lambda n(t))\cos \theta(t)=0$
for all $t \in I$. 
Let $(\overline{\gamma},\overline{\nu}_1,\overline{\nu}_2):I \to \R^3 \times \Delta$ be 
$
\overline{\gamma}(t)=\gamma(t)+\lambda\nu_2(t),\ \overline{\nu}_1(t)=\nu_2(t),\ \overline{\nu}_2(t)=\cos\theta(t)\mu(t)-\sin\theta(t)\nu_1(t).
$
Since $\dot{\overline{\gamma}}(t)=(\alpha(t)+\lambda n(t))\mu(t)-\lambda\ell(t)\nu_1(t)$, we have  $\dot{\overline{\gamma}}(t) \cdot \overline{\nu}_1(t)=\dot{\overline{\gamma}}(t) \cdot \overline{\nu}_2(t)=0$ for all $t \in I$. 
It follows that $(\overline{\gamma},\overline{\nu}_1,\overline{\nu}_2)$ is a framed curve. 
Therefore, $(\gamma,\nu_1,\nu_2)$ is a $(\nu_2,\overline{\nu}_1)$-Bertrand framed curve.
\enD
%%%%%

%%%%%
\begin{proposition}\label{nu2-nu1_curvature}
Suppose that $(\gamma,\nu_1,\nu_2)$ and $(\overline{\gamma},\overline{\nu}_1,\overline{\nu}_2):  I \to \R^3 \times \Delta$ are {$(\nu_2,\overline{\nu}_1)$-mates}, where $\overline{\gamma}(t)=\gamma(t)+\lambda\nu_2(t), \overline{\nu}_2(t)=\cos \theta(t){\mu}(t)-\sin \theta(t){\nu_1}(t), \overline{\mu}(t)=\sin \theta(t){\mu}(t)+\cos \theta(t){\nu_1}(t)$ and $\theta:I \to \R$ is a smooth function. 
Then the curvature $(\overline{\ell}, \overline{m}, \overline{n}, \overline{\alpha})$ of $(\overline{\gamma},\overline{\nu}_1,\overline{\nu}_2)$ is given by
\begin{align*}
  \overline{\ell}(t)&=\ell(t)\sin\theta(t)+n(t)\cos\theta(t), \\
  \overline{m}(t)&=-\ell(t)\cos\theta(t)+n(t)\sin\theta(t), \\
  \overline{n}(t)&=-\dot{\theta}(t)-m(t), \\
  \overline{\alpha}(t)&=(\alpha(t)+\lambda n(t))\sin\theta(t)-\lambda\ell(t)\cos\theta(t). 
\end{align*}
\end{proposition}
%%%%%
\demo
By differentiating $\overline{\nu}_2(t)=\cos \theta(t){\mu}(t)-\sin \theta(t){\nu_1}(t)$, we have 
\begin{align*}
-\overline{\ell}(t)\overline{\nu}_1(t)+\overline{n}(t)\overline{\mu}(t)&=(-\dot{\theta}(t)-m(t))\overline{\mu}(t)+(-n(t)\cos\theta(t)-\ell(t)\sin\theta(t))\nu_2(t).
\end{align*}
Then we have $\overline{n}(t)=-\dot{\theta}(t)-m(t)$. Since $\nu_2(t)=\overline{\nu}_1(t)$, we have $\overline{\ell}(t)=\ell(t)\sin\theta(t)+n(t)\cos\theta(s)$. Moreover, by differentiating $\overline{\mu}(t)=\sin \theta(t){\mu}(t)+\cos \theta(t){\nu_1}(t)$, we have 
\begin{align*}
-\overline{m}(t)\overline{\nu}_1(t)-\overline{n}(t)\overline{\nu}_2(t)&=(\dot{\theta}(t)+m(t))\overline{\nu}_2(t)+(-n(t)\sin\theta(t)+\ell(t)\cos\theta(t))\nu_2(t).
\end{align*}
Since $\nu_2(t)=\overline{\nu}_1(t)$, we have $\overline{m}(t)=-\ell(t)\cos\theta(t)+n(t)\sin\theta(t)$. By $\overline{\alpha}(t)\sin \theta(t)=\alpha(t)+\lambda n(t)$ and $\overline{\alpha}(t)\cos(t)=-\lambda\ell(t)$, we have $\overline{\alpha}(t)=(\alpha(t)+\lambda n(t))\sin\theta(t)-\lambda\ell(t)\cos\theta(t)$.
\enD
%%%%%

%%%%%
\begin{theorem}\label{nu2-nu2}
$(\gamma,\nu_1,\nu_2):I \to \R^3 \times \Delta$ is a $(\nu_2,\overline{\nu}_2)$-Bertrand framed curve if and only if $(\gamma,\nu_1,\nu_2):I \to \R^3 \times \Delta$ is a $(\nu_2,\overline{\nu}_1)$-Bertrand framed curve. 
\end{theorem}
%%%%%
\demo
If $(\gamma,\nu_1,\nu_2):I \to \R^3 \times \Delta$ is a $(\nu_2,\overline{\nu}_2)$-Bertrand framed curve, then there exists a framed curve $(\overline{\gamma},\overline{\nu}_1,\overline{\nu}_2)$ and a smooth function $\lambda:I \to \R$ with $\lambda \not\equiv 0$ such that 
$\overline{\gamma}(t)=\gamma(t)+\lambda(t)\nu_2(t)$ and $\nu_2(t)=\overline{\nu}_2(t)$ for all $t \in I$
Since $(\overline{\gamma},\overline{\nu}_2,\overline{\nu}_1)$ is also a framed curve, $(\gamma,\nu_1,\nu_2)$ is a $(\nu_2,\overline{\nu}_1)$-Bertrand framed curve and vice versa (see, Theorem \ref{relation} and Remark \ref{relation-remark}).
\enD
%%%%%
%%%%%
\begin{remark}{\rm 
Let $(\gamma,\nu_1,\nu_2)$ and $(\overline{\gamma},\overline{\nu}_1,\overline{\nu}_2)$ be $(\nu_2,\overline{\nu}_1)$-mates. 
The curvature of $(\overline{\gamma},\overline{\nu}_1,\overline{\nu}_2)$ is given by $(\overline{\ell}, \overline{m}, \overline{n}, \overline{\alpha})$. 
Then $(\gamma,\nu_1,\nu_2)$ and $(\overline{\gamma},\overline{\nu}_2,\overline{\nu}_1)$ are $(\nu_2,\overline{\nu}_2)$-mates and the curvature of $(\overline{\gamma},\overline{\nu}_2,\overline{\nu}_1)$ is given by $(-\overline{\ell}, -\overline{n}, -\overline{m}, -\overline{\alpha})$ by Propositions \ref{change-frame} and \ref{nu2-nu1_curvature}.
}
\end{remark}
%%%%%

%%%%%
\begin{remark}{\rm 
Let $(\gamma,\nu_1,\nu_2):I \to \R^3 \times \Delta$ is a framed curve with an adapted frame. 
Then curvature of $(\gamma,\nu_1,\nu_2)$ is given by $(0,m,n,\alpha)$ by \eqref{adapted-frame}. 
By Theorems \ref{nu2-nu1} and \ref{nu2-nu2}, $(\gamma,\nu_1,\nu_2)$ is automatically not only $(\nu_2,\overline{\nu}_1)$-Bertrand framed curve, but also $(\nu_2,\overline{\nu}_2)$-Bertrand framed curve.
}
\end{remark}
%%%%%

%%%%%
\begin{theorem}\label{nu2-mu}
$(\gamma,\nu_1,\nu_2):I \to \R^3 \times \Delta$ is a $(\nu_2,\overline{\mu})$-Bertrand framed curve if and only if $\ell(t)=0$ and there exists a smooth function $\lambda : I\rightarrow\R$ with $\lambda \not\equiv 0$ such that $\alpha(t)+\lambda(t)n(t)=0$ for all $t\in I$. 
\end{theorem}
%%%%%
\demo
Suppose that $(\gamma,\nu_1,\nu_2):I \to \R^3 \times \Delta$ is a $(\nu_2,\overline{\mu})$-Bertrand framed curve. 
Then there exists a framed curve $(\overline{\gamma},\overline{\nu}_1,\overline{\nu}_2)$ and a smooth function $\lambda:I \to \R$ with $\lambda \not\equiv 0$ such that 
$\overline{\gamma}(t)=\gamma(t)+\lambda(t)\nu_2(t)$ and $\nu_2(t)=\overline{\mu}(t)$ for all $t \in I$.
By differentiating, we have $\overline{\alpha}(t) \overline{\mu}(t)=(\alpha(t)+\lambda(t)n(t))\mu(t)-\lambda(t)\ell(t)\nu_1(t)+\dot{\lambda}(t)\nu_2(t)=0$ for all $t \in I$. 
Since $\nu_2(t)=\overline{\mu}(t)$, we have $\alpha(t)+\lambda(t)n(t)=0$ and $-\lambda(t)\ell(t)=0$ for all $t \in I$. 
Since $\lambda \not\equiv 0$ and the continuous condition, we have $\ell(t)=0$ for all $t \in I$. 
\par
Conversely, suppose that $\ell(t)=0$ and there exists a smooth function $\lambda : I\rightarrow\R$ such that $\alpha(t)+\lambda(t)n(t)=0$ for all $t\in I$. 
Let $(\overline{\gamma},\overline{\nu}_1,\overline{\nu}_2):I \to \R^3 \times \Delta$ be 
$
\overline{\gamma}(t)=\gamma(t)+\lambda(t)\nu_2(t),\ \overline{\nu}_1(t)=\mu(t),\ \overline{\nu}_2(t)=\nu_1(t).
$
Since $\dot{\overline{\gamma}}(t)=\dot{\lambda}(t)\nu_2(t)$, we have  $\dot{\overline{\gamma}}(t) \cdot \overline{\nu}_1(t)=\dot{\overline{\gamma}}(t) \cdot \overline{\nu}_2(t)=0$ for all $t \in I$. 
It follows that $(\overline{\gamma},\overline{\nu}_1,\overline{\nu}_2)$ is a framed curve. 
Moreover, since $\overline{\mu}(t)=\overline{\nu}_1(t)\times\overline{\nu}_2(t)$, we have $\overline{\mu}(t)={\nu}_2(t)$ for all $t\in I$. 
Therefore, $(\gamma,\nu_1,\nu_2)$ is a $(\nu_2,\overline{\mu})$-Bertrand framed curve.
\enD
%%%%%
%%%%%
\begin{remark} {\rm
If $(\gamma,\nu_1,\nu_2)$ is a $(\nu_2, \overline{\mu})$-Bertrand framed curve and $n(t) \not=0$ for all $t \in I$, then $\overline{\gamma}$ is a circular evolute with respect to $\nu_2$ of the framed curve $(\gamma,\nu_1,\nu_2)$ (cf. \cite{Honda-Takahashi-Preprint}). 
}
\end{remark}
%%%%%
%%%%%
\begin{proposition}\label{nu2-mu_curvature}
Suppose that $(\gamma,\nu_1,\nu_2)$ and $(\overline{\gamma},\overline{\nu}_1,\overline{\nu}_2): I \to \R^3 \times \Delta$ are {$(\nu_2,\overline{\mu})$-mates}, where $\overline{\gamma}(t)=\gamma(t)+\lambda(t)\nu_2(t),  \overline{\nu}_1(t)=\cos \theta(t){\mu}(t)-\sin \theta(t){\nu}_1(t), \overline{\nu}_2(t)=\sin \theta(t){\mu}(t)+\cos \theta(t){\nu}_1(t)$ and $\theta:I \to \R$ is a smooth function. 
Then the curvature $(\overline{\ell}, \overline{m}, \overline{n}, \overline{\alpha})$ of $(\overline{\gamma},\overline{\nu}_1,\overline{\nu}_2)$ is given by
\begin{align*}
\overline{\ell}(t)=-\dot{\theta}(t)-m(t), \ 
\overline{m}(t)=-n(t)\cos\theta(t), \ 
\overline{n}(t)=-n(t)\sin\theta(t), \ 
\overline{\alpha}(t)=\dot{\lambda}(t).
\end{align*}
\end{proposition}
%%%%%
\demo
By differentiating $\overline{\nu}_1(t)=\cos \theta(t){\mu}(t)-\sin \theta(t){\nu}_1(t)$, we have 
\begin{align*}
\overline{\ell}(t)\overline{\nu}_2(t)+\overline{m}(t)\overline{\mu}(t)&=(-\dot{\theta}(t)-m(t))\overline{\nu}_2(t)-n(t)\cos\theta(t)\nu_2(t).
\end{align*}
Then we have $\overline{\ell}(t)=-\dot{\theta}(t)-m(t)$. Since $\nu_2(t)=\overline{\mu}(t)$, we have $\overline{m}(t)=-n(t)\cos\theta(t)$. Moreover, by differentiating $\overline{\nu}_2(t)=\sin \theta(t){\mu}(t)+\cos \theta(t){\nu}_1(t)$, we have 
\begin{align*}
-\overline{\ell}(t)\overline{\nu}_1(t)+\overline{n}(t)\overline{\mu}(t)&=(\dot{\theta}(t)+m(t))\overline{\nu}_1(t)-n(t)\sin\theta(t)\nu_2(t).
\end{align*}
Since $\nu_2(t)=\overline{\mu}(t)$, we have $\overline{n}(t)=-n(t)\sin\theta(t)$. By $\overline{\alpha}(t) \overline{\mu}(t)=(\alpha(t)+\lambda(t)n(t))\mu(t)-\lambda(t)\ell(t)\nu_1(t)+\dot{\lambda}(t)\nu_2(t)=0$ and $\nu_2(t)=\overline{\mu}(t)$, we have $\overline{\alpha}(t)=\dot{\lambda}(t)$.
\enD
%%%%%
\par
Finally, we give concrete examples of Bertrand framed curves.
%%%%%
\begin{example}{\rm
Let $(\gamma,\nu_1,\nu_2):[0,2\pi) \to \R^3 \times \Delta$ be
\begin{align*}
\gamma(t)=\frac{1}{\sqrt{1+p^2}}&\Bigl(-\frac{p}{2}\left(-\frac{1}{q+\sqrt{1+p^2}}\cos(q+\sqrt{1+p^2})t-\frac{1}{q-\sqrt{1+p^2}}\cos(q-\sqrt{1+p^2})t\right),\\
&\quad \frac{p}{2}\left(\frac{1}{q+\sqrt{1+p^2}}\sin(q+\sqrt{1+p^2})t-\frac{1}{q-\sqrt{1+p^2}}\sin(q-\sqrt{1+p^2})t\right),\\
&\quad -\frac{1}{q}\cos qt\Bigr),
\end{align*}
\begin{align*}
\nu_1(t)=\frac{p}{\sqrt{1+p^2}}&\Bigl(-\frac{p}{2}\left(\frac{1}{1+\sqrt{1+p^2}}\sin(1+\sqrt{1+p^2})t+\frac{1}{1-\sqrt{1+p^2}}\sin(1-\sqrt{1+p^2})t\right),\\
&\quad \frac{p}{2}\left(\frac{1}{1+\sqrt{1+p^2}}\cos(1+\sqrt{1+p^2})t-\frac{1}{1-\sqrt{1+p^2}}\cos(1-\sqrt{1+p^2})t\right), \\
&\quad\sin t\Bigr),\\
\nu_2(t)=\frac{p}{\sqrt{1+p^2}}&\Bigl(-\frac{p}{2}\left(-\frac{1}{1+\sqrt{1+p^2}}\cos(1+\sqrt{1+p^2})t-\frac{1}{1-\sqrt{1+p^2}}\cos(1-\sqrt{1+p^2})t\right),\\
&\quad \frac{p}{2}\left(\frac{1}{1+\sqrt{1+p^2}}\sin(1+\sqrt{1+p^2})t-\frac{1}{1-\sqrt{1+p^2}}\sin(1-\sqrt{1+p^2})t\right), \\
&\quad-\cos t\Bigr),
\end{align*}
where $p, q\in \R \setminus \{0\}$ with $q\neq\pm\sqrt{1+p^2}$. 
By a direct calculation,  
\begin{align*}
\mu(t)&=\nu_1(t)\times\nu_2(t)=\frac{1}{\sqrt{1+p^2}}\Bigl(-p\cos\sqrt{1+p^2}t, -p\sin\sqrt{1+p^2}t, 1\Bigr)
\end{align*}
and $(\gamma, \nu_1, \nu_2)$ is a framed curve with the curvature 
$$
(\ell(t), m(t), n(t), \alpha(t))=(0, p\cos t, p\sin t, \sin qt).
$$ 
If we take $\theta(t)=n\pi$ where $n\in \Z$, then the condition of Theorem \ref{framed-Bertrand-equivalent} is satisfied. 
Hence, $(\gamma, \nu_1, \nu_2)$ is a $(\nu_1, \overline{\nu}_1)$ and $(\nu_1, \overline{\nu}_2)$-Bertrand framed curve. 
Similarly, if we take $\theta(t)=\pi/2+n\pi$, the conditions of Theorems \ref{framed-Mannheim-equivalent} and \ref{nu2-nu1} are satisfied. 
Hence, $(\gamma, \nu_1, \nu_2)$ is a $(\nu_2, \overline{\nu}_1)$ and $(\nu_2, \overline{\nu}_2)$-Bertrand framed curve. Moreover, if we take $\theta(t)=\pi/2+n\pi-t$, the conditions of Theorems  \ref{mu-nu1} and \ref{mu-nu2} are satisfied. 
Hence, $(\gamma, \nu_1, \nu_2)$ is a $(\mu, \overline{\nu}_1)$ and $(\mu, \overline{\nu}_2)$-Bertrand framed curve.
\par
If we take $p\not=0, \pm \sqrt{3}$, $q=2$ and $\lambda(t)=-(2/p)\sin t$ (respectively, $\lambda(t)=-(2/p)\cos t$), then the condition of Theorem \ref{nu1-mu} (respectively, Theorem \ref{nu2-mu}) is satisfied. 
Hence, $(\gamma, \nu_1, \nu_2)$ is a $(\nu_1, \overline{\mu})$ (respectively, $(\nu_2, \overline{\mu})$)-Bertrand framed curve.
}
\end{example}
\begin{example}[Spherical Legendre curves \cite{Takahashi}]{\rm
Let $(\gamma,\nu):I \to \Delta$ be a spherical Legendre curve with curvature $(m,n)$.
If we denote $\mu=\gamma \times \nu$, then 
$$
\left(
\begin{array}{c}
\dot{\gamma}(t)\\
\dot{\nu}(t)\\
\dot{\mu}(t)
\end{array} \right)=
\left(
\begin{array}{ccc}
0 & 0 & m(t)\\
0 & 0 & n(t)\\
-m(t) & -n(t) & 0
\end{array}\right)
\left(
\begin{array}{c}
\gamma(t)\\
\nu(t)\\
\mu(t)
\end{array}\right).
$$
It follows that $(\gamma,\gamma,\nu)$ is a framed curve with curvature $(0,m,n,m)$. 
Since $\ell=0$, $(\gamma,\nu_1,\nu_2)=(\gamma,\gamma,\nu)$ is a $(\nu_1,\overline{\nu}_1)$ (respectively, $(\nu_1,\overline{\nu}_2)$, $(\nu_2,\overline{\nu}_1)$, $(\nu_2,\overline{\nu}_2)$)-Bertrand framed curve. 
Moreover, if we take $\lambda=-1$, then the condition of Theorem \ref{nu1-mu} is satisfied. 
Hence, $(\gamma, \nu_1, \nu_2)$ is a $(\nu_1, \overline{\mu})$-Bertrand framed curve. 
\par
It also follows that $(\gamma,\nu,\gamma)$ is a framed curve with curvature $(0,-n,-m,-m)$. 
Since $\ell=0$, $(\gamma,\nu_1,\nu_2)=(\gamma,\nu,\gamma)$ is a $(\nu_1,\overline{\nu}_1)$ (respectively, $(\nu_1,\overline{\nu}_2)$, $(\nu_2,\overline{\nu}_1)$, $(\nu_2,\overline{\nu}_2)$)-Bertrand framed curve. 
Moreover, if we take $\lambda=-1$, then the condition of Theorem \ref{nu2-mu} is satisfied. 
Hence, $(\gamma, \nu_1, \nu_2)$ is a $(\nu_2, \overline{\mu})$-Bertrand framed curve. 
}
\end{example}
%%%%%

%%%%%%%%%% Bibliography %%%%%%%%%%%%%%%%%%

%%%%%%%%
Nozomi Nakatsuyama, 
\\
Muroran Institute of Technology, Muroran 050-8585, Japan,
\\
E-mail address: 23043042@muroran-it.ac.jp
\\
\\
Masatomo Takahashi, 
\\
Muroran Institute of Technology, Muroran 050-8585, Japan,
\\
E-mail address: masatomo@muroran-it.ac.jp

\end{document}